\documentclass[12pt]{amsart}
\usepackage{amssymb}
\usepackage{amsmath}
\usepackage{mathabx}
\usepackage{amsthm}
\usepackage{amstext}
\usepackage{amsopn}
\usepackage{mathrsfs}
\usepackage{latexsym}
\allowdisplaybreaks
\newtheorem{Definition}{Definition}[section]
\newtheorem{Proposition}{Proposition}[section]
\newtheorem{Lemma}{Lemma}[section]
\newtheorem{Theorem}{Theorem}[section]
\newtheorem{Corollary}{Corollary}[section]
\newtheorem{Remark}{Remark}[section]
\newtheorem{Example}{Example}[section]

\begin{document}
\bibliographystyle{plain}
\title[Matricial R-transform]
{Matricial R-transform}
\author[R. Lenczewski]{Romuald Lenczewski}
\address{Romuald Lenczewski, \newline
Instytut Matematyki i Informatyki, Politechnika Wroc\l{}awska, \newline
Wybrze\.{z}e Wyspia\'{n}skiego 27, 50-370 Wroc{\l}aw, Poland  \vspace{10pt}}
\email{Romuald.Lenczewski@pwr.wroc.pl}

\begin{abstract}
We study the addditon problem for strongly matricially free random varia\-bles which 
generalize free random variables.
Using operators of Toeplitz type, we derive a linearization formula for the {\it matricial R-transform} 
related to the associated convolution. It is a linear combination of Voiculescu's 
R-transforms in free probability with coefficients given by internal units of the considered array of subalgebras. 
This allows us to view this formula as the {\it matricial linearization property} of the R-transform. 
Since strong matricial freeness unifies the main types of noncommutative independence, 
the matricial R-transform plays the role of 
a unified noncommutative analog of the logarithm of the Fourier transform for free, boolean, 
monotone, ortho\-gonal, s-free and c-free independence. This paper treats the case of two-dimensional arrays. However,
all results can be generalized to arrays of any finite dimension and a more general version of this paper
will appear elswehere.
\end{abstract}

\maketitle
\section{Introduction}
In this paper we study the addition problem for strongly matricially free random variables [12].
These results extend Voiculescu's results [22] 
on the addition of free random variables (generalized by Maassen [15] and Bercovici and Voiculescu [2])
and include those for the addition of random variables associated 
with other fundamental types of noncommutative independence (monotone [16], boolean [20]), 
generalizations of these (conditionally free [4], conditionally monotone [8])
as well as those related to free subordination (s-free [10], orthogonal [10]).

In classical probability, the addition problem for classically independent random variables 
is related to the classical convolution of probability measures $\mu_1\star\mu_2$ and to the 
lo\-ga\-rithm of the Fourier transform which linearizes this convolution, 
\begin{equation}\tag{1.1}
{\rm log}F_{\mu_1\star\mu_2}={\rm log}F_{\mu_1}+ {\rm log}F_{\mu_2},
\end{equation}
where $F_{\mu}$ is the Fourier transform of $\mu$. However, the analogous addition problem for 
noncommutative random variables is more involved since there is no single notion of 
noncommutative independence.

In free probability [21-25], there is a remarkable analog of the above formula, in which 
the classical convolution is replaced by the free convolution $\mu_1\boxplus \mu_2$
and the role of the logarithm of the Fourier transform is played by the R-transform,
\begin{equation}\tag{1.2}
R_{\mu_1\boxplus\, \mu_2}=R_{\mu_1}+R_{\mu_2},
\end{equation}
where $R_{\mu}$ is the R-transform of $\mu$.
However, the results on additive convolutions associated with other types of noncommutative independence imply
that the R-transform does not retain the linearization property in the general noncommutative framework
and other transforms, such as the reciprocal Cauchy transform [16], the K-transform [20], or the c-free
R-transform [4], have to be used to describe these convolutions.

This situation was one of our motivations to look for 
new types of independence which would unify the existing fundamental types and for
which the associated transform could play the role of a noncommutative logarithm of
the Fourier transform (some motivation came also from [14]). 
In this work, we focus our attention on a generalization of freeness called
{\it strong matricial freeness} [12], in which 
arrays of non-unital *-subalgebras $({\mathcal A}_{i,j})$ of a given 
unital *-algebra ${\mathcal A}$ replace families of unital *-subalgebras
of free probability. We assume, however, that each algebra ${\mathcal A}_{i,j}$ has an internal unit $1_{i,j}$
which is a projection and that all these units, together with $1_{{\mathcal A}}$, 
generate a commutative *-subalgebra ${\mathcal I}$ called the {\it algebra of units}.
Moreover, one distinguished state is replaced by an array of states $(\varphi_{i,j})$
on ${\mathcal A}$, where by a state we understand a complex-valued normalized positive 
linear functional. 

We have shown in [12] that strong matricial freeness unifies 
the fundamental types of noncommutative independence. Namely, including 
some additional cases distinguished in this paper, we obtain a correspondence
$$
{\rm shapes\;of\;arrays}\;\;\rightleftarrows \;\;{\rm types\;of \; independence},
$$
where, in particular, square, lower-triangular and diagonal arrays correspond to freeness, 
monotone independence and boolean independence, respectively. This scheme also includes
the notions of s-freeness and orthogonal independence and leads to the correspondence
$$
{\rm addition\; of \;rows}\;\;\rightleftarrows \;\;{\rm convolutions\; of \;measures},
$$
where by a row we understand the corresponding sum of variables. 
All binary additive convolutions obtained in this fashion are induced by addition of strongly matricially free 
random variables. Since all convolutions considered in this paper will be additive, 
the adjective `additive' will be usually omitted in the sequel.

Our scheme shows that strongly matricially free random variables 
are basic noncommutative variables whose addition leads to well-known 
convolutions. Thus, it is natural to define the 
{\it strongly matricially free convolution} of the array of 
distributions $(\mu_{i,j})$ as the distribution of the sum
\begin{equation}\tag{1.3}
A=\sum_{i,j}a_{i,j}
\end{equation}
of strongly matricially free random variables in a distinguished state $\varphi$, where
$\mu_{i,j}$ is the distribution of $a_{i,j}\in {\mathcal A}_{i,j}$ in the state $\varphi_{i,j}$,
denoted by $\boxplus_{\,i,j}\mu_{i,j}$.
If we work in the category of $C^{*}$-algebras and
the variables are self-adjoint, the considered distributions 
can be identified with compactly supported probability measures on the real line.

For simplicity, let us specialize to the case of arrays $(a_{i,j})$ in which
$(i,j)\in J$, where $J$ is a subset of the Cartesian product $\{1,2\}\times \{1,2\}$, often
omitted in our notations.
The symbol $\boxplus$ used for our convolution operation does not mean that the 
distribution of $A$ is the free convolution of an array of distributions.
In particular, one cannot change the `matricial order' in which the distributions appear. 
However, if the array $(a_{i,j})$ is square and the 
corresponding distributions $(\mu_{i,j})$ are {\it row-identical}, then
\begin{equation}\tag{1.4}
\boxplus_{\,i,j}\mu_{i,j}=\mu_1\boxplus \mu_2,
\end{equation}
where $\mu_{1,1}=\mu_{1,2}=\mu_1$ and $\mu_{2,1}=\mu_{2,2}=\mu_2$. Therefore, in this case,
the new convolution can be viewed as a decomposition of the free convolution. 
Moreover, this framework also gives decompositions of other binary convolutions in terms of $(\mu_{i,j})$
if we consider subarrays of square arrays or relax the assumption that the distributions are row-identical.
In this fashion we obtain boolean, monotone, conditionally free and conditionally monotone convolutions, 
as well as s-free and orthogonal convolutions related to the subordination property of the free convolution 
discovered by Voiculescu [25] and generalized by Biane [3]. For other results on s-free convolutions, 
see also [11] (mulitplicative case) and [17] (multivariate case).

We would like to find a suitable noncommutative analog of the logarithm of the 
Fourier transform. For that purpose we study {\it strongly matricially free Toeplitz ope\-ra\-tors}, 
similar to those in free probability studied by Voiculescu [22] and Haagerup [7] 
(an extension of which to conditional freeness has been proposed recently [18]).
We follow the approach of Haagerup who used adjoints of the Toeplitz operators
introduced in the original approach of Voiculescu. Namely, we take
\begin{equation}\tag{1.5}
a_{i,j}=\ell_{i,j}+f_{i,j}(\ell_{i,j}^{*})
\end{equation}
where $\ell_{i,j}, \ell_{i,j}^{*}$ are strongly matricially free creation and annihilation operators, respectively,
living in the strongly matricially free Fock space and $f_{i,j}$ is a polynomial 
for any $(i,j)\in J$. A novelty, as compared with the free case, is that 
the constant term of $f_{i,j}(\ell_{i,j}^{*})$ is equal to the 
internal unit $1_{i,j}$ multiplied by some complex number.

This allows us to derive an equation for the Cauchy-transform ${\mathcal G}_{A}$ 
of some distribution of $A$ in the state $\varphi$ with 
an operator-valued argument taken from the algebra of units ${\mathcal I}$. 
This equation is of the same form as in the case of scalar-valued arguments, namely
\begin{equation}\tag{1.6}
\mathcal{G}_{A}\left(\frac{1}{z}+\mathcal{R}_{A}(z)\right)=z
\end{equation}
in some `neighborhood' of zero. However, the argument of ${\mathcal G}_{A}$
is an invertible ${\mathcal I}$-valued power series which has a bounded inverse, with ${\mathcal R}_{A}(z)$
playing to role of an ${\mathcal I}$-valued analog of the R-transform. Another important point is that 
the distribution of $A$ considered here, although scalar-valued, cannot be identified with the collection of moments 
of $A$ in the state $\varphi$, but rather with the collection of moments of $A$ alternating with 
elements of ${\mathcal I}$, as in the operator-valued free probability [23,24].

This distribution can be viewed as a noncommutative extension of the strongly matricially free 
convolution considered above, by abuse of notation also denoted $\boxplus_{i,j}\mu_{i,j}$, 
in which we convolve the array $(\mu_{i,j})$ in a more noncommutative fashion.
Using it, we can prove the addition formula for ${\mathcal I}$-valued R-transforms 
involved which not only extends the scalar-valued linearization formula (1.2), 
but also takes a nice form
\begin{equation}\tag{1.7}
\mathcal{R}_{A}=\sum_{i,j}\mathcal{R}_{i,j}
\end{equation}
where ${\mathcal R}_{i,j}(z)=R_{i,j}(z)1_{i,j}$ and $R_{i,j}$ denotes the scalar-valued R-transform of $\mu_{i,j}$ for any $(i,j)\in J$. In fact, this formula shows that convolutions associated with the main types of noncommutative independence can be
linearized by the ${\mathcal I}$-valued R-transform, called {\it matricial R-transform},
except that we have to use noncommutative distributions of $A$ and modify the concept of linearization.
Essentially, this noncommutative setting turns out natural for proving existence of a noncommutative analog 
of (1.1) and (1.2) which is suitable for various types of independence.

We can also view this formula as a {\it linearization property} of the R-transform of a new type,
since the matricial R-transform is a linear combination of the R-transforms of 
Voiculescu with ${\mathcal I}$-valued coefficients. Let us stress that this linear combination
is very special since its coefficients are exactly the internal units of the considered subalgebras. 
We should also add that a similar formula holds for all states from the array $(\varphi_{i,j})$ 
if we take suitably truncated operators and only when considering arrays of Cauchy transforms of this type
we can prove uniqueness of the matricial R-transform.

In particular, if the array is square and distributions are row-identical, 
formula (1.7) gives additivity of the scalar-valued R-transforms since two different 
unit decompositions lead to decompositions of the R-transforms, namely
\begin{equation}\tag{1.8}
1_{\mathcal A}=\sum_{j=1}^{2}1_{i,j}\;\;\;\rightarrow\;\;\;R_{\mu_{i}}=\sum_{j=1}^{2}\mathcal{R}_{i,j}
\end{equation}
for each $i\in \{1,2\}$ and thus (1.2) becomes its consequence. Let us remark in this context that the fact that the distributions are row-identical 
does not mean that such are the states. In fact, in the GNS construction given in [12] 
the off-diagonal states are by construction different from each other and different from the diagonal states.

Conversely, we can arrive at the matricial R-transform starting from free probability
and using operatorial subordination [10]. For that purpose, we begin with a pair of free Toeplitz operators
$a_1,a_2\in B({\mathcal F}({\mathcal H}))$, where ${\mathcal F}({\mathcal H})$ is the full Fock space over a two-dimensional Hilbert space with an orthonormal basis $\{e_1,e_2\}$. The results of [10] give decompositions 
\begin{equation}\tag{1.9}
a_{j}=x_j+X_j \;\;\;{\rm for}\;\;j\in \{1,2\}
\end{equation}
where 
$x_j$ and $X_j$ are equal to $a_j$ restricted to ${\mathcal F}({\mathbb C}e_{j})$ 
and its orthogonal complement ${\mathcal F}({\mathbb C}e_{j})^{\perp}$, respectively, where $j\in \{1,2\}$. 
In this case, the $\varphi$-distribution of $x_1$ agrees with the $\varphi_2$-distribution 
of $X_1$ and the $\varphi$-distribution of $x_2$ agrees with the $\varphi_1$-distribution of 
$X_2$ (these distributions form a matrix and can be denoted $\mu_{1,1}$, $\mu_{1,2}$, $\mu_{2,2}$ and $\mu_{2,1}$, respectively), where $\varphi_j$ is the state defined by the vector $e_j$ and $j\in \{1,2\}$. 

Treating the R-transforms of the distributions of $a_1$ and $a_2$ as operators 
on ${\mathcal F}({\mathcal H})$ and decomposing the unit of $B({\mathcal F}({\mathcal H}))$
as described above, we arrive at (1.7) with $A=a_1+a_2$. Then it remains to relax the assumption that 
the distributions are row-identical to obtain the general case.
Therefore, another way to arrive at the matricial R-transform is to follow the 
`subordination path'
$$
\left(
\begin{array}{l}
a_1\\
a_2
\end{array}
\right)
\;\;\;\rightarrow \;\;\;
\left(
\begin{array}{ll}
x_1 & X_1\\
X_2 & x_2
\end{array}
\right)\;\;\;
\rightarrow 
\;\;\;
\left(
\begin{array}{ll}
{\mathcal R}_{1,1} & {\mathcal R}_{1,2}\\
{\mathcal R}_{2,1} & {\mathcal R}_{2,2}
\end{array}
\right),
$$
namely first decompose free random variables as in (1.9) and then assign to them 
the associated operatorial R-transforms induced by unit decompositions (1.8). 

If we deform these R-transforms and thus allow for
different R-transforms in each row, the variables corresponding to rows become
conditionally free with respect to $(\varphi, \psi)$, where $\psi$ is any extension of $\varphi_1$ 
and $\varphi_2$. In this context, let us point out that the advantage of using strong matricial freeness 
instead of conditional freeness is that the first one includes monotone independence, whereas
the second one requires additional (rather restrictive) assumptions on the considered subalgebras, 
although satisfied in the case of algebras of polynomials used when studying convolutions [6].

At the same time, our theory has some new features since it provides a framework placed
between scalar- and operator-valued free probability.
An important property is that despite the fact that the internal units are not identified, 
the considered states are scalar-valued. A similar approach was used to define {\it matricial freeness} [12], 
which could also be called `weak matricial freeness', especially when compared with strong matricial freeness.
Addition of (weakly) matricially free random variables will be studied elsewhere.

Let us remark that our construction of a noncommutative logarithm unifying the logarithm of the Fourier transform and the K-transform given in [9] is of a different nature than the approach presented in this work, although it is possible that 
our present approach can be extended in that direction.

The remainder of this paper is organized as follows. Section 2 is devoted to the concept of strong matricial freeness.
The corresponding convolution is studied in Section 3. Suitable Toeplitz operators living in the strongly 
matricially free Fock space are introduced and studied in Section 4. In Section 5, we define 
the matricial R-transform and prove the corresponding linearization formula, 
which is the main result of this paper. In Section 6, we give conditions under 
which the matricial R-transform is unique. In Sections 7 and 8, we study 
strongly matricially free convolutions from a combinatorial point of view.
\section{Strong matricial freeness}

Let us recall the basic notions related to the concept of strong matricial freeness [12].
Let ${\mathcal A}$ be a unital *-algebra with an array $({\mathcal A}_{i,j})$ 
of non-unital *-subalgebras of ${\mathcal A}$ and let $(\varphi_{i,j})$ be a 
family of states on ${\mathcal A}$ (there is some similarity to freeness with infinitely many states [5]).
Further, we assume that each ${\mathcal A}_{i,j}$ has an {\it internal unit} $1_{i,j}$ which is a projection 
for which $a1_{i,j}=1_{i,j}a=a$ for any $a\in {\mathcal A}_{i,j}$ and
that the unital subalgebra ${\mathcal I}$ of ${\mathcal A}$ generated by all internal 
units is commutative.

Crucial is the fact that internal units are not identified with the unit of ${\mathcal A}$ and therefore additional 
conditions on moments involving these units are needed.
In order to state these conditions, we 
need subsets of $(I\times I)^{m}$ of the form
\begin{eqnarray*}
\Gamma_{m}&=&\{((i_1,i_2), (i_2,i_3), \ldots , (i_m,i_{m+1})): i_1\neq i_2\neq \ldots \neq i_m\}
\end{eqnarray*}
where $I$ is an index set and $m\in {\mathbb N}$, with the corresponding union denoted 
$$
\Gamma=\bigcup_{m=1}^{\infty}\Gamma_m.
$$
Objects (algebras and their elements, states, etc.) labelled by indices of the form $(j,j)$ 
for any $j$ and $i\neq j$ for any $i\neq j$, respectively, will be called {\it diagonal} and
{\it off-diagonal}.
Important is the difference between diagonal and off-diagonal objects. 
\begin{Definition}
{\rm 
We say that $(1_{i,j})$ is a {\it strongly matricially free array of units} associated with 
$({\mathcal A}_{i,j})$ and $(\varphi_{i,j})$ if for any diagonal state $\varphi$ it holds that
\begin{enumerate}
\item[(a)]
$\varphi(u_1au_2)=\varphi(u_1)\varphi(a)\varphi(u_2)$ 
for any $a\in {\mathcal A}$ and $u_1,u_2\in {\mathcal I}$,
\item[(b)]
$\varphi(1_{i,j})=\delta_{i,j}$ for any $i,j$,
\item[(c)]
if $a_{k}\in {\mathcal A}_{i_k,j_k}\cap {\rm Ker}\varphi_{i_k,j_k}$, where $1<k\leq m$, then
$$
\varphi(a1_{i_1,j_1}a_{2}\ldots a_m)=
\left\{
\begin{array}{cc}
\varphi(aa_{2} \ldots a_n) & {\rm if}\;\;((i_{1},j_{1}), \ldots , (i_{m},j_m))\in \Gamma\\
0 & {\rm otherwise}
\end{array}
\right..
$$
where $a\in {\mathcal A}$ is arbitrary and $(i_1,j_1)\neq \ldots \neq (i_m,j_m)$.
\end{enumerate}}
\end{Definition}

The main condition of strong matricial freeness reminds freeness as 
the definition gives below shows.

\begin{Definition}
{\rm 
We say that *-subalgebras $({\mathcal A}_{i,j})$ are {\it strongly} 
{\it matricially free} with respect to $(\varphi_{i,j})$ if the array of internal units 
$(1_{i,j})$ is the associated strongly matricially free array of units and
\begin{equation}\tag{2.1}
\varphi(a_{1}a_{2}\ldots a_{n})=0\;\;\; {\rm whenever }\;\;a_{k}\in {\mathcal A}_{i_k,j_k}\cap {\rm Ker}\varphi_{i_k,j_k}
\end{equation}
for any diagonal state $\varphi$, where $(i_1,j_1)\neq \ldots \neq (i_n,j_n)$.}
\end{Definition}

Using these conditions, one can easily show that multiplication of weak and strong matricially free random variables in the kernel form reminds multiplication of matrices and arrays of units remind matrix units, except that in the strong case the elements which are on their diagonals survive only in the last matrix.
Naturally, the array of variables $(a_{i,j})$ in a unital *-algebra ${\mathcal A}$ is called 
strongly matricially free with respect to $(\varphi_{i,j})$ if
there exists an array of projections $(1_{i,j})$ which is a 
strongly matricially free array of units associated with ${\mathcal A}$ and $(\varphi_{i,j})$ 
and such that the array $({\mathcal A}_{i,j})$ of *-subalgebras 
of the form ${\mathcal A}_{i,j}=alg(a_{i,j},a_{i,j}^{*},1_{i,j})$ 
is strongly matricially free with respect to $(\varphi_{i,j})$.
If ${\mathcal A}$ is a $C^{*}$-algebra, we assume that the subalgebras ${\mathcal A}_{i,j}$ and 
${\mathcal I}$ are $C^{*}$-subalgebras of ${\mathcal A}$.

We shall consider an array of states $(\varphi_{i,j})$ on ${\mathcal A}$ in which diagonal states coincide with a distinguished state $\varphi$, whereas the off-diagonal ones are given by $\varphi_{i,j}=\varphi_j$, where $i\neq j$ and 
\begin{equation}\tag{2.2}
\varphi_{j}(a)=\varphi(b_{j}^{*}ab_{j})
\end{equation}
for some $b_{j}\in {\mathcal A}_{j,j}\cap {\rm Ker}\varphi$ such that 
$\varphi(b_{j}^{*}b_{j})=1$, any $a\in {\mathcal A}$ and $j$, called
{\it conjugate states} (they were called `conditions' in [12]). 
In particular, this implies the normalization conditions
$\varphi_{j}(1_{i,k})=\delta_{j,k}$ for any $i,j,k$.

Finally, if the distribution of $a_{i,j}$ in the state $\varphi_{i,j}$ does not depend 
on $j$, we will say that the array $(a_{i,j})$ is {\it row-identically distributed}.

\begin{Proposition}
Let $({\mathcal A}_{i,j})$ be strongly matricially free with respect to $(\varphi_{i,j})$, where
$\varphi_{j,j}=\varphi$ for any $j$ and $\varphi_{i,j}=\varphi_j$ for any $i\neq j$.  
Let us consider the mixed moment $\varphi_{i,j}(a_{1}a_{2}\ldots a_{n})$, 
where $a_{k}\in {\mathcal A}_{i_k,j_k}$ and $(i_1,j_1)\neq \ldots \neq (i_n,j_n)$. 
\begin{enumerate}
\item
If $a_n$ is diagonal, $j_n\neq j$ and the state $\varphi_{i,j}$ is off-diagonal, then the moment vanishes.
In particular, the off-diagonal state $\varphi_{i,j}$ vanishes on any ${\mathcal A}_{k,k}$ 
for $k\neq j$.
\item
If $a_n$ is off-diagonal and the state $\varphi_{i,j}$ is diagonal, then the moment vanishes. 
In particular, a diagonal state vanishes on any off-diagonal subalgebra.
\item
If $a_k$ is diagonal for some $1\leq k <n$ and $a_{r}\in {\rm Ker}\varphi_{i_r,j_r}$ 
for $k<r\leq n$, then the moment vanishes.
\end{enumerate}
\end{Proposition}
{\it Proof.}
These are straightforward consequences of the definition of strong matricial freeness.
\hfill $\blacksquare$
\begin{Definition}
{\rm 
By the {\it strongly matricially free Fock space} over 
the array $({\mathcal H}_{i,j})$
we understand the Hilbert space of the form
\begin{equation}\tag{2.3}
\mathcal{N}=
{\mathbb C}\Omega \oplus \bigoplus _{m=1}^{\infty}
\bigoplus_{\stackrel {i_{1}\neq \ldots \neq i_{m}}
{\scriptscriptstyle n_1, \ldots , n_m \in {\mathbb N}}}
{\mathcal H}_{i_1,i_2}^{\otimes \,n_1}
\otimes 
{\mathcal H}_{i_2,i_3}^{\otimes \,n_2}
\otimes \ldots \otimes 
{\mathcal H}_{i_{m},i_{m}}^{\otimes \,n_m},
\end{equation}
where $\Omega$ is a unit vector, with the canonical inner product.}
\end{Definition}

In particular, when the array $({\mathcal H}_{i,j})$ is square and consists
of one-dimensional Hilbert spaces ${\mathcal H}_{i,j}={\mathbb C}e_{i,j}$ for $i,j\in \{1,2\}$, we have
$$
\mathcal{N}
=\bigoplus_{m=0}^{\infty}\mathcal{N}^{(m)},
$$
where the first few summands are of the form
\begin{eqnarray*}
\mathcal{N}^{(0)}&=&{\mathbb C}\Omega\\
\mathcal{N}^{(1)}&=&{\mathbb C}e_{1,1}
\oplus {\mathbb C}e_{2,2}\\
\mathcal{N}^{(2)} &=&
{\mathbb C}e_{1,1}^{\otimes 2}\oplus 
{\mathbb C}e_{2,2}^{\otimes 2}\oplus
{\mathbb C}(e_{1,2}\otimes e_{2,2})\oplus
{\mathbb C}(e_{2,1}\otimes e_{1,1}) \\
\mathcal{N}^{(3)}&=&
{\mathbb C}e_{1,1}^{\otimes 3}
\oplus 
{\mathbb C}e_{2,2}^{\otimes 3}
\oplus 
{\mathbb C}(e_{2,1}\otimes e_{1,1}^{\otimes 2})
\oplus 
{\mathbb C}(e_{1,2}\otimes e_{2,2}^{\otimes 2})
\oplus 
{\mathbb C}(e_{2,1}^{\otimes 2}\otimes e_{1,1})\\
&&
\oplus\,
{\mathbb C}(e_{1,2}^{\otimes 2}\otimes e_{2,2})
\oplus 
{\mathbb C}(e_{1,2}\otimes e_{2,1}\otimes e_{1,1})
\oplus 
{\mathbb C}(e_{2,1}\otimes e_{1,2}\otimes e_{2,2}),
\end{eqnarray*}
etc. Note that, in contrast to the (weakly) matricially free Fock space,
diagonal vectors appear only at the ends of the tensor products. 

Let us distinguish the closed subspaces ${\mathcal N}_{i,j}$ spanned by simple tensors which 
begin with $e_{i,j}$, where $(i,j)\in J$. Then we have a decomposition
\begin{equation}\tag{2.4}
{\mathcal N}={\mathbb C}\Omega \oplus \bigoplus_{i,j}{\mathcal N}_{i,j}
\end{equation}
and the diagonal unit $1_{j,j}$ is the projection onto 
${\mathbb C}\Omega \oplus {\mathcal N}_{j,j}\cong {\mathcal F}({\mathbb C}e_{j,j})$, 
whereas the off-diagonal unit $1_{i,j}$ is the projection onto ${\mathcal N}\ominus {\mathcal F}({\mathbb C}e_{i,i})$.
\begin{Remark}
{\rm Let $p$ be the projection onto ${\mathbb C}\Omega$ and let 
$p_{i,j}$ be a projection onto ${\mathcal N}_{i,j}$ for any $i,j$.
Note that we have `orthogonal decompositions' of the unit
from ${\mathcal A}=B({\mathcal N})$, 
\begin{equation}\tag{2.5}
1_{{\mathcal A}}=1_{1,1}+1_{1,2}=1_{2,2}+1_{2,1},
\end{equation}
as well as
\begin{equation}\tag{2.6}
1_{{\mathcal A}}=p +\sum_{i,j}p_{i,j}=\sum_{i,j}q_{i,j},
\end{equation}
where $q_{1,1}=p$, $q_{1,2}=p_{2,2}$, $q_{2,1}=p_{1,1}$ and $q_{2,2}=p_{1,2}+p_{2,1}$,
which will be of use in our study of the matricial R-transform.}
\end{Remark}
\begin{Definition}
{\rm Let $A=(\alpha_{i,j})$ be an array of positive real numbers
and let $({\mathcal H}_{i,j})=({\mathbb C}e_{i,j})$ be the associated array of
Hilbert spaces. By the {\it strongly matricially free creation operators} associated with $A$ we understand
operators of the form
\begin{equation}\tag{2.7}
\ell_{i,j}=\alpha_{i,j}\sigma^{*}\ell(e_{i,j})\sigma, 
\end{equation}
where $\sigma: {\mathcal N}\rightarrow {\mathcal F}(\bigoplus_{i,j}{\mathcal H}_{i,j})$
is the canonical embedding in the given free Fock space and
the $\ell(e_{i,j})$'s denote the canonical free creation operators.
By the {\it strongly matricially free annihilation operators} we understand their adjoints.}
\end{Definition}

\begin{Proposition}
The array $({\mathcal A}_{i,j})$ of *-subalgebras of $B({\mathcal N})$, 
where ${\mathcal A}_{i,j}=alg(\ell_{i,j}, \ell_{i,j}^{*})$ for any $(i,j)\in J$, is strongly matricially 
free with respect to $(\varphi_{i,j})$.
\end{Proposition}
{\it Proof.}
The proof is left to the reader since it is very similar to that in the free case (see also [13, Proposition 4.1]).
Let us only observe that the relation 
\begin{equation}\tag{2.8}
\ell_{i,j}^{*}\ell_{i,j}^{}=\alpha_{i,j}^{2}1_{i,j}^{}
\end{equation}
implies that 
$1_{i,j}\in {\mathcal A}_{i,j}$ since $\alpha_{i,j}$ is assumed to be positive.
\hfill  $\blacksquare$

\begin{Remark}
{\rm 
In a similar way one defines the concept of {\it matricial freeness}, which, 
especially when compared with the strong case, can also be called `weak matricial freeness'. 
The difference is that in the definition of the array of matricially free array of units we replace the sets
$\Gamma_{n}$ by slightly larger sets
$$
\Lambda_{m}=\{((i_1,j_1), (i_2,j_2), \ldots , (i_m,j_{m})): (i_1, i_2)\neq (i_2,i_3)\neq \ldots \neq (i_m,i_{m+1})\},
$$
where all `matricially related' pairs of indices are allowed. In this case, in order to obtain a non-zero
mixed moment of type  $\varphi(a_1a_2\ldots a_n)$ it is also necessary that $a_n$ is diagonal, but 
diagonal variables can also appear in the middle of non-zero moments of the type considered in Proposition 2.1.}
\end{Remark}

\section{Convolutions}

From the results of [12] it follows that different shapes of $(a_{i,j})$ correspond 
to different types of noncommutative independence. The convolution which unifies the associated
convolutions is the convolution of an array of distributions defined below.

\begin{Definition}
{\rm Let $(a_{i,j})$ be a two-dimensional array of strongly matricially free random variables
with the corresponding array of distributions $(\mu_{i,j})$ in the states $(\varphi_{i,j})$. 
The distribution of the sum
\begin{equation}\tag{3.1}
A=\sum_{i,j}a_{i,j}
\end{equation}
in the state $\varphi$ is called the {\it strongly matricially free convolution} and 
will be denoted $\boxplus_{\,i,j}\mu_{i,j}$ by analogy with the free convolution.
In these notations it is understood that the summation runs over pairs $(i,j)\in J$, 
which is often omitted.}
\end{Definition}

In particular, the strongly matricially free convolution generalizes the free convolution, which
justifies our notation. For simplicity, let us consider a two-dimensional array $(a_{i,j})$ 
of random variables and an array of states given by 
\begin{equation}\tag{3.2}
\left(
\begin{array}{cc} 
\varphi_{1,1}   & \varphi_{1,2}\\
\varphi_{2,1} & \varphi_{2,2}
\end{array}
\right) 
=
\left(
\begin{array}{ll} 
\varphi   & \varphi_2\\
\varphi_1 & \varphi
\end{array}
\right) 
\end{equation}
where $\varphi$ is a distinguished state on ${\mathcal A}$ and $\varphi_1,\varphi_2$ are conjugate states, 
with respect to which the considered array is strongly matricially free.

For the sake of greater generality, by a two-dimensional array we 
shall understand a subbarray of a square two-dimensional array. 
Consider sums 
\begin{equation}\tag{3.3}
a_{i}=\sum_{j}a_{i,j}
\end{equation}
where $i=1,2$ and it is understood that the summation runs over those $j$'s, for which $(i,j)\in J$. 
In other words, we add variables which appear in each row separately.
Clearly, these sums may reduce to one term if only one variable appears in a given row.
Let us remark that the array is not assumed to contain the diagonal as 
in the case of limit theorems [12,13]. This allows
us to consider arrays corresponding to s-free independence and orthogonal independence.

In the theorem given below we state explicitly how addition of rows corresponds to
convolutions associated with various notions of noncommutative independence 
and how these convolutions can be decomposed using the strongly matricially free convolution.
By the $\varphi$-distributions of $A$ and $a_{i}$ we understand the collections 
of moments of $A$ and $a_{i}$, respectively, in the state $\varphi$.
Each pair of variables is denoted $\{a_1,a_2\}$, but one should remember that in 
some cases (3-5) order in which they appear is relevant. 
\begin{Theorem}
Assume that $(a_{i,j})$ is strongly matricially free with respect to 
$(\varphi_{i,j})$ and row-identically distributed with distributions $(\mu_{i,j})$. Let
$\mu$ and $\mu_i$ be the $\varphi$-distributions of $A$ and $a_i$, where $i\in \{1,2\}$, 
respectively.
\begin{enumerate}
\item
If the array is square, then $\{a_{1},a_{2}\}$ is free with respect to $\varphi$
and $\mu=\mu_1 \boxplus \mu_2$.
\item
If the array is diagonal, then $\{a_1,a_2\}$ is boolean independent with respect to $\varphi$ and
$\mu=\mu_1\uplus \mu_2$.
\item
If the array is lower-triangular, then $\{a_{1},a_{2}\}$ is monotone independent with respect
to $\varphi$ and $\mu=\mu_1\vartriangleright \mu_2$.
\item
If the array is upper-anti-triangular, then $\{a_1,a_2\}$ is s-free independent with respect to $(\varphi, \varphi_1)$
and $\mu=\mu_1 \boxright \mu_2$.
\item
If the array consist of one column, then $\{a_1,a_2\}$ is orthogonally independent with respect 
to $(\varphi, \varphi_1)$ and $\mu=\mu_1\vdash \mu_2$.
\end{enumerate}
\end{Theorem}
{\it Proof.} Since $A=a_{1}+a_{2}$, cases (1)-(3) are immediate consequences
of [12, Proposition 4.1], where it was shown that $\{a_1,a_2\}$ is free, boolean independent 
or monotone independent, depending on whether the array is square, diagonal or lower-triangular, respectively.
In order to show cases (4) and (5) it suffices to show that $\{a_{1,1}+a_{1,2}, a_{2,1}\}$ 
is s-free and $\{a_{1,1}, a_{2,1}\}$
is orthogonal under $(\varphi, \varphi_1)$. However, these properties follow from the GNS representation 
in the subordination context given in [10].
\hfill $\blacksquare$\\

We have distinguished only those convolutions which are non-trivial and cannot be reduced to the other ones. 
This eliminates upper-triangular arrays (formally they correspond to anti-monotone 
independence and its conditional counterpart), which can be reduced to lower-triangular ones,
as well as {\it lower-anti-triangular} arrays (i.e. those in which there is no $a_{1,1}$), 
which can be reduced to {\it upper-anti-triangular} ones (i.e those in which there is no $a_{2,2}$).
Finally, we have omitted arrays consisting of one row, in which case the $\varphi$-distribution of 
$A$ agrees with the $\varphi$-distribution of the diagonal variable $a_{1,1}$ and 
the $\varphi_2$-distribution of $A$ agrees with that of the off-diagonal variable $a_{1,2}$ 
(all mixed moments of two variables vanish).

If we lift the assumption that the variables are row-identically distributed, we obtain 
the conditionally free additive convolution.
In a similar manner, a generalization of 
the case of lower-triangular arrays leads to conditionally monotone variables of Hasebe [8].

\begin{Theorem} 
If $(a_{i,j})$ is a square array of random variables in ${\mathcal A}$ which is
strongly matricially free with respect to $(\varphi_{i,j})$ and has distributions $(\mu_{i,j})$, then 
\begin{equation}\tag{3.4}
\boxplus_{\,i,j}\mu_{i,j}=(\mu_{1,1}, \mu_{1,2})\boxplus (\mu_{2,2}, \mu_{2,1})
\end{equation}
and $\{a_1,a_2\}$ is conditionally free with respect to $(\varphi, \psi)$, where
$\psi$ is any state on ${\mathcal A}$ extending $\varphi_1|_{{\mathcal A}_{2}}$ and 
$\varphi_2|_{{\mathcal A}_{1}}$, where ${\mathcal A}_{1}=alg({\mathcal A}_{1,1},{\mathcal A}_{1,2})$ 
and ${\mathcal A}_{2}=alg({\mathcal A}_{2,1}, {\mathcal A}_{2,2})$.
\end{Theorem}
{\it Proof.}
It suffices to prove that variables of the form
$$
a_{i}^{0}=a_{i}-\psi(a_{i})=
\sum_{j}(a_{i,j}-\psi(a_{i,j}))
$$
satisfy the equation of c-freeness with respect to $\varphi$, namely
$$
\varphi(a_{i_1}^{0}a_{i_2}^{0}\ldots a_{i_n}^{0})=
\varphi(a_{i_1}^{0})\varphi(a_{i_2}^{0})\ldots \varphi(a_{i_n}^{0})
$$
where $i_1\neq i_2\neq \ldots \neq i_n$ and $\psi$ is any extension of 
$\varphi_1|_{{\mathcal A}_{2}}$ and $\varphi_2|_{{\mathcal A}_{1}}$. 
The assumption on $\psi$ gives
\begin{eqnarray*}
\varphi(a_{1}^{0})&=&\sum_{j}(\varphi(a_{1,j})-\psi(a_{1,j}))=\varphi(a_{1,1})-\psi(a_{1,2})\\
\varphi(a_{2}^{0})&=&\sum_{j}(\varphi(a_{2,j})-\psi(a_{2,j}))=\varphi(a_{2,2})-\psi(a_{2,1}),
\end{eqnarray*}
since $\varphi(a_{1,2})=\varphi(a_{2,1})=0$ as well as $\varphi_2(a_{1,1})=\varphi_1(a_{2,2})=0$
by Proposition 2.1.
Therefore, our variables can be decomposed as
\begin{eqnarray*}
a_{i}^{0}&=&a_{i}'+a_{i}''+b_{i}
\end{eqnarray*}
where the first two terms belong to diagonal subalgebras, namely
$$
a_{i}'=a_{i,i}-\varphi(a_{i,i})1_{i,i}
$$
for each $i$ and
$$
a_{1}''=(\varphi(a_{1,1})-\psi(a_{1,2}))1_{1,1},\;\;\;
a_{2}''=(\varphi(a_{2,2})-\psi(a_{2,1}))1_{2,2}
$$
whereas the third one is an element of an off-diagonal one, i.e.
$$
b_1=a_{1,2}-\psi(a_{1,2})1_{1,2} \;\;\;{\rm and}\;\;\; b_2=a_{2,1}-\psi(a_{2,1})1_{2,1}.
$$
Now, using Proposition 2.1(2) again, we obtain
$$
\varphi(a_{i_1}^{0}a_{i_2}^{0}\ldots a_{i_n}^{0})=
\varphi(a_{i_1}^{0}a_{i_2}^{0}\ldots a_{i_n}')+
\varphi(a_{i_1}^{0}a_{i_2}^{0}\ldots a_{i_n}'').
$$
Next, we claim that 
$$
\varphi(a_{i_1}^{0}a_{i_2}^{0}\ldots a_{i_n}')=0.
$$
Namely, in view of Proposition 2.1(3), the diagonal variables can appear 
only as the last ones in non-zero mixed moments of kernel form. 
This implies, however, that 
$$
\varphi(a_{i_1}^{0}a_{i_2}^{0}\ldots a_{i_n}')=
\varphi(b_{i_1}b_{i_2}\ldots b_{i_{n-1}}a_{i_n}')=0
$$
by strong matricial freeness (note that this kind of reduction is not possible 
in the case of weak matricial freeness). It remains to use property (a) of 
Definition 2.1 to write
$$
\varphi(a_{i_1}^{0}a_{i_2}^{0}\ldots a_{i_n}'')=
(\varphi(a_{i_n,i_n})-\psi(a_{i_n,i_n}))\varphi(a_{i_1}^{0}a_{i_2}^{0}\ldots a_{i_{n-1}}^{0})
$$
and use an inductive argument to finish the proof.
\hfill $\blacksquare$

\section{Toeplitz operators}
In order to study the addition of free random variables, Voiculescu used Toeplitz operators [22].
They turn out to be closely related to R-transforms and can be used to prove that 
they are additive under free convolutions. A new proof of additivity was presented by Haagerup [7] who 
used the adjoints of Toeplitz operators of the form
\begin{equation}\tag{4.1}
a=\ell_{1} + f(\ell_{1}^{*}) \;\;\;{\rm and}\;\;\;b=\ell_{2}+g(\ell_{2}^{*})
\end{equation}
where $\ell_1, \ell_2$ are isometries on the full Fock space ${\mathcal F}({\mathcal H})$, where
${\mathcal H}$ is a two-dimensional Hilbert space with orthonormal basis $\{e_1,e_2\}$,
given by tensoring on the left by $e_1$ and $e_2$, respectively, and where $f,g$ are 
polynomials. Using only Fock space techniques, he derived the addition formula for the R-transforms. 

We would like to generalize this approach to the case of strong matricial freeness in order to find a suitable analog
of the R-transform.
\begin{Definition}
{\rm Let $(\ell_{i,j})$ be the array of strongly matricially free creation operators on ${\mathcal N}$ 
and let $f_{i,j}$ be a polynomial for any $(i,j)\in J$. Operators of the form
\begin{equation}\tag{4.2}
a_{i,j}=\ell_{i,j}+f_{i,j}(\ell_{i,j}^{*})
\end{equation}
where $(i,j)\in J$ and the constant term of $f_{i,j}$ is the internal unit $1_{i,j}$
multiplied by a complex number, will be called {\it strongly matricially free Toeplitz operators}.}
\end{Definition}
These operators are similar to free Toeplitz operators 
used by Voiculescu and Haagerup. We will need the distributions of these operators in the states
$(\varphi_{i,j})$ defined by the array of unit vectors $(\Omega_{i,j})$, where 
$$
\Omega_{j,j}=\Omega\;\;\;{\rm and}\;\;\; \Omega_{i,j}=e_{j,j}\;\;{\rm for}\;\;i\neq j,
$$
which replace the single vacuum vector which suffices in the free case.
It is easy to see that the restriction of $a_{i,j}$ to the closed subspace
of ${\mathcal N}$ spanned by vectors 
$$
\{\ell_{i,j}^{n}\Omega_{i,j}:n\geq 0\}
$$
has the Toeplitz form for any $(i,j)\in J$. Moreover, in the case of row-identically distributed square arrays
these operators are closely related to those in free probability since the sum of those lying 
in the $j$-th row give a decomposition of the $j$-th free Toeplitz operator (see Introduction).
\begin{Proposition}
The R-transform of the distribution $\mu_{i,j}$ of the operator
$a_{i,j}$ in the state $\varphi_{i,j}$ is given by
\begin{equation}\tag{4.3}
R_{i,j}(z)=f_{i,j}(\alpha_{i,j}^{2}z),
\end{equation}
where $(i,j)\in J$ and the constant term of $f_{i,j}$ is a complex number.
\end{Proposition}
{\it Proof.}
The proof is similar to that in [7], but some modifications need to be 
introduced. For fixed $(i,j)\in J$, take the vector 
$$
\rho_{i,j}(z)=(1-z\ell_{i,j})^{-1}\Omega_{i,j}=\Omega_{i,j}+\sum_{n=1}^{\infty}z^{n}\alpha_{i,j}^{n}e_{i,j}^{\otimes n}\otimes \Omega_{i,j}
$$
where we set $e_{j,j}^{\otimes n}\otimes \Omega\equiv e_{j,j}^{\otimes n}$ for any $j$ and assume that
$|\alpha_{i,j}z|<1$ in order to get a convergent series. Then
$$
\ell_{i,j}\rho_{i,j}(z)=\frac{1}{z}(\rho_{i,j}(z)-\Omega_{i,j})
$$
for $0<|z|<|\alpha_{i,j}|^{-1}$.
Moreover, since $\ell_{i,j}^{*}\Omega_{i,j}=0$ and 
$$
\ell_{i,j}^{*}(e_{i,j}^{\otimes n}\otimes \Omega_{i,j})=\alpha_{i,j}(e_{i,j}^{\otimes (n-1)}\otimes \Omega_{i,j}),
$$
we can see that 
$$
\ell_{i,j}^{*}\rho_{i,j}(z)=z\alpha_{i,j}^{2}\rho_{i,j}(z)
$$
and thus
$$
a_{i,j}\rho_{i,j}(z)=\frac{1}{z}(\rho_{i,j}(z)-\Omega_{i,j})+f_{i,j}(\alpha_{i,j}^{2}z)\rho_{i,j}(z),
$$
which leads to the equation
$$
\left(\frac{1}{z}+f_{i,j}(\alpha_{i,j}^{2}z)-a_{i,j}\right)\rho_{i,j}(z)=\frac{1}{z}\Omega.
$$
Now, if $|z|$ is sufficiently small and positive, the operator on the left hand side is invertible, which gives
$$
G_{i,j}\left(\frac{1}{z}+f_{i,j}(\alpha_{i,j}^{2}z)\right)=z,
$$
where $G_{i,j}$ denotes the Cauchy transform of the $\varphi_{i,j}$-distribution of  
$a_{i,j}$, and that implies that its R-transform has the desired form by the uniqueness of the R-transform.
\hfill $\blacksquare$\\

It should be noted that in the study of $\varphi$-distributions of the strongly matricially free convolution
it suffices to use strongly matricially free Toeplitz operators $(a_{i,j})$ and their sum $A$ since, in view of 
Proposition 2.2, the array $(a_{i,j})$ is strongly matricially free with respect 
to $(\varphi_{i,j})$ and this implies [12, Propositions 2.2-2.3] that the mixed moments of these variables are uniquely determined by the marginal distributions and thus by the corresponding R-transforms $(R_{i,j})$. 
Therefore, if we start with an arbitrary array of strongly 
matricially free random variables and the same distributions as 
those of $(a_{i,j})$, 
then we obtain the same mixed moments. Thus, without loss of generality, we can 
use strongly matricially free Toeplitz operators with given R-transforms. 

It can also be justifed that in this study one can use polynomials as R-transforms instead of infinite series
$$
R_{i,j}(z)=\sum_{n=1}^{\infty}r_{i,j}(n)z^{n-1}
$$
for any $(i,j)\in J$. This is because any mixed moment of order $m$ of the sum of arbitrary strongly matricially free
random variables can be expressed in terms of a finite number of coefficients $r_{i,j}(n)$ for $1\leq n \leq m$.

These results give us a hint that one might be able to use strongly matricially free Toeplitz operators
to obtain an analog of the R-transform. Following Haagerup's approach [7], 
we shall first define suitable vectors of geometric series type and then act on these vectors
with Toeplitz operators. 

Let $(\ell_{i,j})$ be the array of strongly matricially free creation operators 
living in the strongly matricially free Fock space ${\mathcal N}$ and let
$$
L=\sum_{i,j}\ell_{i,j}
$$
be the corresponding sum of creation operators. Consider the vector
$$
\rho(z)=(1-zL)^{-1}\Omega,
$$
where $z$ is a complex number such that $|z|<(\sum_{i,j}|\alpha_{i,j}|^2)^{-1}$, which ensures convergence of the associated series.

\begin{Lemma}
The sum $A=\sum_{i,j}a_{i,j}$ of strongly matricially free Toeplitz operators satisfies the equation
\begin{equation}\tag{4.4}
A\rho(z)=\frac{1}{z}(\rho(z)-\Omega)+\sum_{i,j}f_{i,j}(\alpha_{i,j}^{2}z)1_{i,j}\rho(z)
\end{equation}
where $\rho(z)=(1-zL)^{-1}\Omega\;$ and $\;0<|z|<(\sum_{i,j}|\alpha_{i,j}|^{2})^{-1}$.
\end{Lemma} 
{\it Proof.}
Clearly, as in the free case, we have
$$
L\rho(z)=\frac{1}{z}(\rho(z)-\Omega).
$$
Let us now decompose the vector $\rho(z)$ into sums of two vectors in two different ways. 
These vectors will be eigenvectors of annihilation operators from each row of the array $(\ell_{i,j}^{*})$. 
Thus, using the same notation as in the proof of Proposition 4.1, we get the action of the diagonal annihilation operators
$$
\ell_{j,j}^{*}\rho_{j,j}^{}(z)=\alpha_{j,j}^{2}z\rho_{j,j}^{}(z)
$$
for $j\in \{1,2\}$. Moreover, for the off-diagonal ones we have
$$
\ell_{i,j}^{*}\eta_{i,j}(z)=\alpha_{i,j}^{2}z\eta_{i,j}(z),
$$
where $i\neq j$ and 
\begin{eqnarray*}
\eta_{1,2}(z)&=&\rho(z)-\rho_{1,1}(z)\\
\eta_{2,1}(z)&=&\rho(z)-\rho_{2,2}(z).
\end{eqnarray*}
Let us show the first of these equations. 
We have
\begin{eqnarray*}
\ell_{1,2}^{*}(\rho(z)-\rho_{1,1}(z))&=& \sum_{n=1}^{\infty}z^{n}(\ell_{1,2}^{*}L^{n}-\ell_{1,2}^{*}\ell_{1,1}^{n})\Omega\\
&=&\sum_{n=1}^{\infty}z^{n}\ell_{1,2}^{*}\ell_{1,2}^{}L^{n-1}\Omega\\
&=&\alpha_{1,2}^{2}z\sum_{n=1}^{\infty}z^{n-1}(L^{n-1}\Omega -e_{1,1}^{\otimes (n-1)})\\
&=&\alpha_{1,2}^{2}z(\rho(z)-\rho_{1,1}(z))
\end{eqnarray*}
since 
$\ell_{1,2}^{*}\ell_{1,2}=\alpha_{1,2}^{2}1_{1,2}$ and $1_{1,2}L^{n-1}\Omega=L^{n-1}\Omega - e_{1,1}^{\otimes (n-1)}$,
which proves the desired equation.
Finally, 
$$
\ell_{1,2}^{*}\rho_{1,1}^{}(z)=\ell_{2,1}^{*}\rho_{2,2}^{}(z)=
\ell_{1,1}^{*}\eta_{1,2}(z)=\ell_{2,2}^{*}\eta_{2,1}(z)=0.
$$
Therefore,
\begin{eqnarray*}
(f_{1,1}^{}(\ell_{1,1}^{*})+f_{1,2}(\ell_{1,2}^{*}))\rho(z)&=& f_{1,1}(\alpha_{1,1}^{2}z)\rho_{1,1}(z)+
f_{1,2}(\alpha_{1,2}^{2}z)\eta_{1,2}(z)\\
(f_{2,2}^{}(\ell_{2,2}^{*})+f_{2,1}(\ell_{2,1}^{*}))\rho(z)&=& f_{2,2}(\alpha_{2,2}^{2}z)\rho_{2,2}(z)+
f_{2,1}(\alpha_{2,1}^{2}z)\eta_{2,1}(z),
\end{eqnarray*}
which gives (4.4) and thus completes the proof.
\hfill $\blacksquare$
\begin{Lemma}
Let $\varphi$ be the state associated with the vacuum vector $\Omega$. Then there exists $\epsilon$ such that
\begin{equation}\tag{4.5}
z=\varphi\left(\left(\frac{1}{z}+\sum_{i,j}R_{i,j}(z)1_{i,j}-A\right)^{-1}\right)
\end{equation}
whenever $0<|z|<\epsilon$.
\end{Lemma}
{\it Proof.}
Rewriting the result of Lemma 4.1, we obtain
$$
\left(\frac{1}{z}+\sum_{i,j}f_{i,j}(\alpha_{i,j}^{2}z)1_{i,j}-A\right)\rho(z)=\frac{1}{z}\Omega.
$$
For small and positive $|z|$, the operator on the LHS is invertible and then 
$$
z\rho(z)=\left(\frac{1}{z}+\sum_{i,j}f_{i,j}(\alpha_{i,j}^{2}z)1_{i,j}-A\right)^{-1}\Omega
$$
which, with the use of $\langle z\rho(z), \Omega \rangle =z$ and Proposition 4.1, finishes the proof. 
\hfill $\blacksquare$\\

The above lemma indicates that a natural analog of the R-transform 
could have a more noncommutative character. This is really the case and, 
as we argue in the next section, in order to obtain the corresponding linearization formula, we need
to use more noncommutative distributions of $A$ than those studied in Section 3. The latter 
are then obtained by restriction to the algebra of polynomials in $A$.

\section{Matricial R-transform}

In view of Lemma 4.2 and its similarity to the scalar-valued counterpart involving the scalar-valued R-transform, 
we are ready to define an object called an `operatorial R-transform'. For that purpose we shall use 
the notion of a distribution of $A$ in the state $\varphi$ similar to that in the operator-valued free probability [23].

Let $a$ be an element of a unital complex Banach algebra ${\mathcal A}$ with a 
unital closed subalgebra ${\mathcal I}$ and let $\varphi$ be a normalized linear functional on ${\mathcal A}$. 
The algebra freely generated by ${\mathcal I}$ and an indeterminate $X$ will be denoted by 
${\mathcal I}\langle X\rangle $. The linear functional $\mu:{\mathcal I}\langle X \rangle\rightarrow {\mathbb C}$ 
defined by $\mu=\varphi\circ \tau$, where $\tau:{\mathcal I}\langle X \rangle\rightarrow {\mathcal A}$ 
is the unique homomorphism such that $\tau(b)=b$ for $b\in {\mathcal I}$ and $\tau(X)=a$
will be called the {\it distribution} of $a\in {\mathcal A}$ in the state $\varphi$.
Quantities of the form 
$$
\varphi(b_{n_1}ab_{n_2}\ldots b_{n_{m-1}}ab_{n_m}), {\rm where}\;b_{n_k}\in {\mathcal I}\;\;{\rm for}\;\; 1\leq k \leq m\; 
{\rm and} \;m\in {\mathbb N}\
$$
will be called the ${\mathcal I}$-{\it moments} of $a$. In particular, the collection of moments 
$$
\{\varphi(b(ab)^{n}): n\in {\mathbb N}\cup \{0\}\}
$$ 
will be called the {\it $b$-distribution of }$a$. Note that we do not assume, as in the operator-valued free probability [23], where $\varphi$ is a conditional expectation, that $\varphi(b)=b$ for any $b\in {\mathcal I}$.

It is clear that if $b\in {\mathcal I}$ is invertible and $\parallel b^{-1}\parallel <\parallel a \parallel^{-1}$ 
then the inverse of $b-a$ exists and takes the form of a series
$$
(b-a)^{-1}=\sum_{n=0}^{\infty}b^{-1}(ab^{-1})^{n}
$$
which converges in the norm topology. This leads to the operatorial analog 
of the Cauchy transform of the form 
\begin{equation}\tag{5.1}
\mathcal{G}_{a}(b)=\sum_{n=0}^{\infty}\varphi\left(b^{-1}(ab^{-1})^{n}\right),
\end{equation}
due to continuity of $\varphi$, which plays the role of the Cauchy transform of the $b$-distribution of $A$ 
in the state $\varphi$. If ${\mathcal A}$ is a $C^{*}$-algebra, continuity follows from the fact that $\varphi$ is a state.
Clearly, in this case $b, b^{-1}$ and $a$ become bounded operators on some Hilbert space ${\mathcal H}$.

\begin{Definition}
{\rm Under the above assumptions, let $\mu$ denote the distribution of $a\in {\mathcal A}$
in the state $\varphi$.
If there exists an ${\mathcal I}$-valued power series of the form
\begin{equation}\tag{5.2}
\mathcal{R}_{a}(z)=\sum_{n=1}^{\infty}c_{n}z^{n-1},
\end{equation}
where $c_{n}\in {\mathcal I}$ for all $n\in {\mathbb N}$ and $z\in {\mathbb C}$, which is 
convergent in the norm topology for sufficiently small $|z|$, and for which it holds that
\begin{equation}\tag{5.3}
\mathcal{G}_{a}\left(\frac{1}{z}+\mathcal{R}_{a}(z)\right)=z
\end{equation}
whenever $|z|$ is sufficiently small and positive, it will be called 
an {\it operatorial R-transform of the distribution $\mu$}.}
\end{Definition}
\begin{Remark}
{\rm Let ${\mathcal A}$ be the $C^{*}$-subalgebra generated by the Toeplitz operator $a_{i,j}$ and 
the unit $1_{i,j}$ and let ${\mathcal I}={\mathbb C}[1_{i,j}]$, where the pair $(i,j)\in J$ is fixed. 
Then
$$
{\mathcal G}_{a_{i,j}}\left(\frac{1}{z}+\mathcal{R}_{a_{i,j}}(z)\right)=z
$$
for small and positive $|z|$, where ${\mathcal R}_{{a_{i,j}}}(z)=R_{i,j}(z)1_{i,j}$ and $R_{i,j}$ is given by Proposition 4.1,
thus ${\mathcal R}_{a_{i,j}}$ becomes an operatorial R-transform of the $\varphi_{i,j}$-distribution of $a_{i,j}$.
}
\end{Remark}

On the level of formal power (Laurent) series, invertibility of the argument of the Cauchy transform 
of the form (5.3) is shown as in the case of scalar-valued arguments. 
Moreover, if $z{\mathcal R}_{A}(z)$ is a contraction, 
the inverse converges in the norm topology. These facts are  
elementary and we present them without a proof.
\begin{Proposition}
Formal power series $C(z)=\sum_{n=0}^{\infty}c_{n}z^{n-1}$ and $B(z)=\sum_{n=0}^{\infty}b_{n}z^{n+1}$, where
$c_{n},b_{n}\in {\mathcal I}$ for any $n$ and $c_{0}=b_{0}=1$, are multiplicative inverses if and only if
$$
b_{m}=\sum_{p=1}^{m}(-1)^{p}\sum_{n_1+\ldots +n_p=m}c_{n_1}\ldots c_{n_p}
$$
$$
c_{m}=\sum_{p=1}^{m}(-1)^{p}\sum_{n_1+\ldots + n_p=m}b_{n_1}\ldots b_{n_p}
$$
where $m,n_1, \ldots , n_p\in {\mathbb N}$. If $\parallel zC(z)-1\parallel <1$, then the power series $B(z)$ 
converges in the norm topology. If $\parallel B(z)-z\parallel <1$, then the series $C(z)$ converges in 
the norm topology.
\end{Proposition}

Note that if ${\mathcal R}_{a}$ is an operatorial R-transform associated with ${\mathcal G}_{a}$, then the operatorial argument of ${\mathcal G}_{a}$ is of the form
\begin{equation}\tag{5.4}
C_{a}(z)=\frac{1}{z}+{\mathcal R}_{a}(z)
\end{equation}
and plays the role analogous to that of the right composition inverse of the Cauchy transform. However, 
since ${\mathcal R}_{a}$ is operator-valued, it is not unique. 

\begin{Definition}
{\rm Let $A\in {\mathcal A}$ be the sum of random variables $(a_{i,j})$ in a unital complex $C^{*}$-algebra ${\mathcal A}$
which are strongly matricially free with respect to $(\varphi_{i,j})$ 
and let ${\mathcal I}$ be its unital $C^{*}$-subalgebra generated  by the internal units. If an 
${\mathcal I}$-valued operatorial R-transform $\mathcal{R}_{A}$ of the $\varphi$-distribution of $A$ takes the form
\begin{equation}\tag{5.5}
{\mathcal R}_{A}(z)=\sum_{i,j}{\mathcal R}_{i,j}(z)
\end{equation}
where ${\mathcal R}_{i,j}(z)=R_{i,j}(z)1_{i,j}$ and $R_{i,j}$ is the R-transform of $\mu_{i,j}$, the distribution of
$a_{i,j}$ in the state $\varphi_{i,j}$, it will be called a {\it matricial R-transform} of the $\varphi$-distribution 
of $A$.}
\end{Definition}

If $\mathcal{R}_{A}$ is a matricial R-transform, then the multiplicative inverse of $C_{A}$ can be
expressed in terms of R-transforms of free convolutions of $(\mu_{i,j})$, with  
orthogonal projections $(q_{i,j})$ introduced in Remark 2.1 being the coefficients.
In particular, this is the case when $(a_{i,j})$ is an array of strongly matricially
free Toeplitz operators, but the result is more general. 
\begin{Proposition}
If $(a_{i,j})$ is the array of Toeplitz operators (4.2) and ${\mathcal R}_{A}$
is the corresponding matricial R-transform, then the multiplicative inverse of 
$C_{A}(z)$ takes the form
$$
B_{A}(z)=\sum_{i,j}\left(\frac{z}{1+zQ_{i,j}(z)}\right)q_{i,j}
$$
where $|z|$ is sufficiently small and positive and $(Q_{i,j})$ is the 
array of R-transforms of free convolutions of the distributions of $(a_{i,j})$ in the states 
$(\varphi_{i,j})$.
\end{Proposition}
{\it Proof.}
Decomposing the internal units $(1_{i,j})$ in terms of $(q_{i,j})$ 
defined in Remark 2.1, we arrive at the orthogonal decomposition
$$
C_{A}(z)=\sum_{i,j}\left(\frac{1}{z}+Q_{i,j}(z))\right)q_{i,j},
$$
where $(Q_{i,j})$ is the array of the form
$$
\left(
\begin{array}{ll}
Q_{1,1}& Q_{1,2}\\
Q_{2,1}& Q_{2,2}
\end{array}
\right)
=
\left(
\begin{array}{ll}
R_{1,1}+R_{2,2}& R_{1,1}+R_{2,1}\\
R_{2,2}+R_{1,2}& R_{1,2}+R_{2,1}
\end{array}
\right),
$$
and thus, by additivity of the R-transform, it is the array of R-transforms 
of free convolutions of $(\mu_{i,j})$, the distributions of $(a_{i,j})$ in the states
$(\varphi_{i,j})$, respectively, namely
$$
\left(
\begin{array}{ll}
\mu_{1,1}\boxplus \mu_{2,2}& \mu_{1,1}\boxplus \mu_{2,1}\\
\mu_{2,2}\boxplus \mu_{1,2}& \mu_{1,2}\boxplus \mu_{2,1}
\end{array}
\right).
$$
This allows us to find the explicit form of the inverse of $C_{A}(z)\in B({\mathcal N})$, namely
$$
B_{A}(z)=\sum_{i,j}\left(\frac{z}{1+zQ_{i,j}(z)}\right)q_{i,j}
$$
which converges in the norm topology for sufficiently small and positive $|z|$ since each $Q_{i,j}$ is analytic in some neighborhood 
of zero, as claimed.
\hfill $\blacksquare$\\

\begin{Remark}
{\rm It is not hard to show that Proposition 5.2 can be generalized to the case when 
$(a_{i,j})$ is an arbitrary array of strongly matricially free variables from 
a unital complex $C^{*}$-algebra ${\mathcal A}$. 
In that case, the projections $(q_{i,j})$ are to be understood 
as suitable elements of ${\mathcal I}$ expressed in terms of $(1_{i,j})$ by the same formulas as in the case
of ${\mathcal N}$.}
\end{Remark}
\begin{Lemma}
An ${\mathcal I}$-valued power series $\mathcal{R}(z)=\sum_{n=1}^{\infty}c_{n}z^{n-1}$ 
converging in the norm topology in a neighborhood of zero is 
an operatorial R-transform of the $\varphi$-distribution of $A$ if and only if
\begin{equation}\tag{5.6}
\sum_{k=1}^{m}\sum_{n_1+\ldots +n_k=m-k}
\varphi(b_{n_1}Ab_{n_2}\ldots b_{n_{k-1}}Ab_{n_k})=0
\end{equation}
for all $m\geq 2$, where we assume that $n_1, \ldots , n_k$ are non-negative integers and where the series $B(z)=\sum_{n=0}^{\infty}b_{n}z^{n+1}$ is the multiplicative inverse of 
$C(z)=1/z+\mathcal{R}(z)$.
\end{Lemma}
{\it Proof.}
In (5.1) we substitute the explicit form for the inverse of $b=C(z)$
given by the series $B(z)$ of Proposition 5.1. Note that if $|z|$ is small, then 
the condition $\parallel zC(z)-1\parallel =\parallel z\mathcal{R}(z)\parallel <1$ is satisfied.
Moreover, if $|z|$ is sufficiently small, the norm of $B(z)$ can be made smaller than $\parallel A\parallel^{-1}$, 
which gives absolute convergence of the series (5.1) in some neighborhood of zero.
By continuity of $\varphi$ this leads to the identity
\begin{eqnarray*}
z
&=&\sum_{n_1=0}^{\infty}\varphi(b_{n_1})z^{n_1+1}+\sum_{n_1,n_2=0}^{\infty}\varphi(b_{n_1}Ab_{n_2})z^{n_1+n_2+2}\\
&+&\sum_{n_1,n_2,n_3=0}^{\infty}\varphi(b_{n_1}Ab_{n_2}Ab_{n_3})z^{n_1+n_2+n_3+3} + \ldots\\
&=&\sum_{m=1}^{\infty}\sum_{k=1}^{m}\sum_{n_1+\ldots +n_k=m-k}\varphi(b_{n_1}Ab_{n_2}\ldots b_{n_{k-1}}Ab_{n_k})z^{m}
\end{eqnarray*}
for small $|z|$, which gives necessary and sufficient conditions for ${\mathcal R}(z)$ to be an 
operatorial R-transform of the $\varphi$-distribution of $A$, namely 
$$
\sum_{k=1}^{m}
\sum_{n_1+\ldots +n_k=m-k}
\varphi(b_{n_1}Ab_{n_2}\ldots b_{n_{k-1}}Ab_{n_k})=0
$$
for all natural $m$, where we assume that the summation indices $n_1, \ldots , n_p$ are 
non-negative integers, which completes the proof.
\hfill $\blacksquare$\\

\begin{Remark}
{\rm 
It is easy to see that (5.6) is a recursion since it is equivalent to 
$$
\varphi(b_{m-1})=-\sum_{k=2}^{m}
\sum_{n_1+\ldots +n_k=m-k}
\varphi(b_{n_1}Ab_{n_2}\ldots b_{n_{k-1}}Ab_{n_k})
$$
for $m\geq 2$ and $b_{0}=1$, and it is clear that it always has a solution.
The conditions for the first few coefficients $b_{m}$ take the form
\begin{eqnarray*}
\varphi(b_1)&=&-\varphi(A)\\
\varphi(b_2)&=&-\varphi(Ab_1)-\varphi(b_1A)-\varphi(A^2)\\
\varphi(b_3)&=&-\varphi(Ab_2)-\varphi(b_1Ab_1)-\varphi(b_2A)-\varphi(A^2b_1)-\varphi(b_1A^2)-\varphi(Ab_1A)-\varphi(A^3)
\end{eqnarray*}
etc., with $b_{0}=1$. In the scalar-valued case, we obtain a unique solution for all coefficients
$b_{m}$ which are expressed in terms of $\varphi$-moments of $A$. In the operator-valued case, conjugate states 
will be used to address the uniqueness problem in Section 6.}
\end{Remark}

We are ready to state and prove the main result of this paper, the addition formula for 
the matricial R-transform, which can also be viewed as the {\it matricial linearization property of the R-transform}.
\begin{Theorem}
If $(a_{i,j})$ is an array of random variables from a unital complex $C^{*}$-algebra ${\mathcal A}$
which is strongly matricially free with respect to $(\varphi_{i,j})$
and $(R_{i,j})$ is the corresponding array of R-transforms, then
\begin{equation}\tag{5.7}
\mathcal{R}_{A}(z)=\sum_{i,j}\mathcal{R}_{i,j}(z),
\end{equation}
where $A=\sum_{i,j}a_{i,j}$ and $\mathcal{R}_{i,j}(z)=R_{i,j}(z)1_{i,j}$ for any $(i,j)\in J$, with 
sufficiently small $|z|$, is an operatorial R-transform of the $\varphi$-distribution of $A$.
\end{Theorem}
{\it Proof.}
If $(a_{i,j})$ is the array of Toeplitz operators, the assertion is a consequence of Lemma 4.2.
It remains to be shown that the case of arbitrary arrays $(a_{i,j})$ reduces to that of Toeplitz operators.
Suppose $(a_{i,j})$ is an array of elements of a $C^{*}$-algebra ${\mathcal A}$ which is strongly matricially free under 
$(\varphi_{i,j})$ with $(1_{i,j})$ being the associated array of internal units and let $(R_{i,j})$ be
the array of R-transforms of the corresponding distributions.
Let $R_{i,j}$ be analytic for $|z|<\epsilon_{i,j}$, where $(i,j)\in J$ and each $\epsilon_{i,j}$ is a positive number. Then the series
$$
{\mathcal R}_{A}(z)=\sum_{i,j}R_{i,j}(z)1_{i,j}=\sum_{n=1}^{\infty}\left(\sum_{i,j}r_{i,j}(n)1_{i,j}\right)z^{n-1}
$$
converges in the norm topology to an element of ${\mathcal I}$ for $|z|<\epsilon$, where 
$\epsilon=\min_{i,j}\epsilon_{i,j}$. For each natural $p$, consider its truncation
$$
{\mathcal R}_{A^{(p)}}(z)=\sum_{i,j}{\mathcal R}_{i,j}^{(p)}(z),
$$
where 
$$
{\mathcal R}_{i,j}^{(p)}(z)=R_{i,j}^{(p)}(z)1_{i,j}\;\;\;{\rm and}\;\;\;
R_{i,j}^{(p)}(z)=\sum_{n=1}^{p}r_{i,j}(n)z^{n-1}
$$
for any $(i,j)\in J$ and $|z|<\epsilon$, and where $A^{(p)}=\sum_{i,j}a_{i,j}^{(p)}$ is the sum of
strongly matricially free Toeplitz operators with $R$-transforms $(R_{i,j}^{(p)})$, respectively. 
By Lemma 4.2, we have
$$
{\mathcal G}_{A^{(p)}}\left(\frac{1}{z}+{\mathcal R}_{A^{(p)}}(z)\right) =z
$$
for any $p$ and $0<|z|<\delta\leq \epsilon$ and some $\delta$.
Therefore, by Lemma 5.1, 
$$
\sum_{k=1}^{m}
\sum_{n_1+\ldots +n_k=m-k}
\varphi(b_{n_1}^{(p)}A^{(p)}b_{n_2}^{(p)}\ldots b_{n_{k-1}}^{(p)}A^{(p)}b_{n_k}^{(p)})=0
$$
for all $m\geq 2$, with $b_{0}^{(p)}=1$. 
Now, by Proposition 5.2, $b_{n}^{(p)}$'s are the coefficients of the power series
$$
B^{(p)}(z)=\sum_{i,j}\left(\frac{z}{1+zQ_{i,j}^{(p)}(z)}\right)q_{i,j}
$$
where $Q_{i,j}^{(p)}$ is the truncation of $Q_{i,j}$ to the polynomial of order $p-1$ and $(q_{i,j})$
is the array of projections of Remark 2.1.
On the other hand, $b_{n}$'s are the coefficients of the power series
$$
B(z)=\sum_{i,j}\left(\frac{z}{1+zQ_{i,j}(z)}\right)q_{i,j},
$$
where, by abuse of notation, we use the same symbols for $(q_{i,j})$, 
which are now the projections from ${\mathcal I}$ (expressed in terms of $(1_{i,j})$ by suitable equations, 
the same as in the case of ${\mathcal N}$).
It can be seen from Proposition 5.1 or directly from the above formulas that $b_n^{(p)}=b_n^{}$ whenever $0\leq n \leq p$.
Moreover, similar arguments as those concerning the moments of $A$ 
in the state $\varphi$ can be used for the alternating mixed moments of $A$ (since $b_n$'s and $b_{n}^{(p)}$'s are linear combinations of internal units)
to conclude that 
$$
\varphi(b_{n_1}^{(p)}A^{(p)}b_{n_2}^{(p)}\ldots b_{n_{k-1}}^{(p)}A^{(p)}b_{n_k}^{(p)})
=
\varphi(b_{n_1}Ab_{n_2}\ldots b_{n_{k-1}}Ab_{n_k})
$$
whenever $n_1+\ldots + n_k+k-1\leq p$. In this fashion, condition (5.6) for any given $m$ can be shown
to be satisfied by taking sufficiently large $p$. This implies that
the series ${\mathcal R}_{A}$ is an operatorial R-transform associated with ${\mathcal G}_{A}$,
which completes the proof.
\hfill $\blacksquare$\\
 
In contrast to highly non-unique operatorial R-transforms of Definition 5.1, 
the matricial R-transform can be made unique if we consider the whole array of states $(\varphi_{i,j})$ 
and not only the distinguished state $\varphi$.
Before we discuss the uniqueness problem, we shall make a few simple remarks 
related to the cases discussed in Theorems 3.1-3.2.
\begin{enumerate}
\item
If the array is square and row-identically distributed, then the matricial R-transform 
associated with ${\mathcal G}_{A}$ takes the form
$$
{\mathcal R}_{A}(z)=
R_{\mu_1}(z)1_{{\mathcal A}} + R_{\mu_2}(z)1_{{\mathcal A}}
$$
due to unit decompositions (1.8) and can be identfied
with the scalar-valued R-transform of $\mu_1\boxplus \mu_2$.
\item
If the array is diagonal, then the matricial 
R-transform associated with ${\mathcal G}_{A}$ takes the form
$$
{\mathcal R}_{A}(z)=
R_{\mu_1}(z)1_{1,1}+ R_{\mu_2}(z)1_{2,2}
$$
which linearizes the extended boolean convolution.
\item
If the array is lower-triangular and row-identically distributed, then the matricial 
R-transform associated with ${\mathcal G}_{A}$ takes the form
$$
{\mathcal R}_{A}(z)=
R_{\mu_1}(z)1_{1,1}+ R_{\mu_2}(z)1_{{\mathcal A}}
$$
which linearizes the extended monotone convolution.
\item
In the general case and with the use of notations used for c-free convolutions, we can write 
the matricial R-transform associated with ${\mathcal G}_{A}$ as
$$
\mathcal{R}_{A}(z)=\sum_{i,j}{\mathcal Q}_{i,j}(z),
$$ 
where ${\mathcal Q}_{i,j}(z)=Q_{i,j}(z)q_{i,j}$ for any $i,j$ and $(Q_{i,j})$ 
is the array of R-transforms of free convolutions
$$
\left(
\begin{array}{ll}
\mu_{1}\boxplus \mu_{2} & \mu_{1}\boxplus \nu_{2}\\
\mu_{2}\boxplus \nu_{1} & \nu_{1}\boxplus \nu_{2}
\end{array}
\right),
$$
obtained by rewriting the array given in the proof of Proposition 5.2.
\end{enumerate}

Without any further assumptions, the matricial R-transform in the state $\varphi$ is not unique. 
In most cases this fact is easy to see since the R-transform of $A$, say 
$R_{A}$, always exists and the corresponding scalar-valued $B(z)$ 
satisfies conditions (5.6), but it is usually different 
than ${\mathcal R}_{A}$ given by (5.7).
However, if the array is square and distributions are row-identical, then 
$\mathcal{R}_{A}$ agrees with $R_{A}$ since $\sum_{j}1_{i,j}=1_{{\mathcal A}}$ for any $i\in \{1,2\}$ 
(then, to show non-uniqueness, one needs to find another example).

\section{Uniqueness}
Let us investigate the uniqueness problem for the operatorial R-transform. For that purpose we shall 
study higher order states, including the conjugate states. 

Anticipating that some essential work will have to be done for operators of Toeplitz type, we first
introduce vectors in the Fock space ${\mathcal N}$ which are similar to $\omega(z)$, but which are generated by 
vectors $e_{1,1}$ and $e_{2,2}$. Let
\begin{equation}\tag{6.1}
\omega_{j}(z)=(1-zL)^{-1}e_{j,j}
\end{equation}
for any $j\in \{1,2\}$, where $|z|$ is sufficiently small to ensure convergence.

\begin{Proposition}
Let $f$ be a polynomial and let the constant term of $f(\ell^{*}_{i,j})$ be a constant 
multiplied by $1_{i,j}$ for any $(i,j)\in J$. Then
\begin{equation}\tag{6.2}
f(\ell^{*}_{j,j})\omega_{j}^{}(z)=f(\alpha_{j,j}^{2}z)1_{j,j}^{}\omega_{j,j}^{}(z)+\alpha_{j,j}Df^{}(\alpha_{j,j}^{2}z)\Omega
\end{equation}
for any $j\in \{1,2\}$, where $Df(z)=\Delta f(z)/z$, with $\Delta f(z)=f(z)-f(0)$. Moreover,
\begin{equation}\tag{6.3}
f(\ell^{*}_{i,j})\omega_{k}(z)=f_{i,j}(\alpha_{i,j}^{2}z)1_{i,j}\omega_{k}(z)
\end{equation}
whenever $i\neq j$ and $k$ is arbitrary.
\end{Proposition}
{\it Proof.}
The proof is similar to that of Lemma 4.1. For instance,
\begin{eqnarray*}
\ell_{j,j}^{*}\omega_{j}(z)&=&
\ell_{j,j}^{*}e_{j,j}+z\ell_{j,j}^{*}\ell_{j,j}^{}(1+zL+z^2L^{2}+\ldots )e_{j,j}\\
&=&\alpha_{j,j}\Omega+z\alpha_{j,j}^{2}1_{j,j}^{}\omega_j^{}(z)
\end{eqnarray*}
for any $j\in \{1,2\}$, which leads to the recursion
$$
(\ell_{j,j}^{*})^{n}\omega_{j}^{}(z)=\alpha_{j,j}^{}(z\alpha_{j,j}^{2})^{n-1}\Omega+
(z\alpha_{j,j}^{2})^{n}1_{j,j}^{}\omega_{j}^{}(z)
$$
for any natural $n$ since $\ell_{j,j}^{*}$ commutes with $1_{j,j}$. This, in turn, gives
$$
f(\ell^{*}_{j,j})\omega_{j}^{}(z)=f(\alpha_{j,j}^{2}z)1_{j,j}^{}\omega_{j}^{}(z)+\alpha_{j,j}Df^{}(\alpha_{j,j}^{2}z)\Omega
$$
for any $j$ and any polynomial $f$, which proves (6.2). The proof of (6.3) is exactly the same as that of Lemma 4.1. 
\hfill $\blacksquare$\\

In addition to the Cauchy transform studied in Section 5, 
we shall now study the Cauchy transforms of noncommutative distributions of $a\in {\mathcal A}$ in the conjugate 
states $\varphi_1$ and $\varphi_2$. This leads to the array
\begin{equation}\tag{6.4}
\mathcal{G}_{a}^{(i,j)}\left(b\right)=\varphi_{i,j}((b-a)^{-1})=\sum_{n=0}^{\infty}\varphi_{i,j}(b^{-1}(ab^{-1})^{n})
\end{equation}
whenever $(i,j)\in J$ and $b\in {\mathcal I}$ is invertible with $\parallel b^{-1} \parallel < \parallel a \parallel^{-1}$.
In particular, we would like to study the case when $a$ is of the form
\begin{equation}\tag{6.5}
A_{i,j}:=P_{i,j}AP_{i,j}, \;\;{\rm where}\;\;(i,j)\in J
\end{equation}
and $A=\sum_{i,j}a_{i,j}$, with $(P_{i,j})$ being the array of projections given by
\begin{equation}\tag{6.6}
P_{i,j}=
\left\{
\begin{array}{lc}
1_{{\mathcal A}}-1_{1,1}1_{2,2}& {\rm if} \;\;i=j
\\
1_{{\mathcal A}}-1_{j,j}& {\rm if}\;\;i \neq j
\end{array}
\right..
\end{equation}
Since these projections belong to ${\mathcal I}$, the operators $A_{i,j}$ belong to ${\mathcal A}$ for any $(i,j)\in J$. 
A more intuitive understanding of these notions can be acquired in the Fock-space context, to which we return below.

\begin{Lemma}
Let $A$ be the sum of strongly matricially free Toeplitz operators in ${\mathcal N}$ with the matricial
R-transform ${\mathcal R}_{A}$ given by (5.7). Then
\begin{equation}\tag{6.10}
\mathcal{G}_{i,j}\left(\frac{1}{z}+\mathcal{R}_{A}(z)\right)=z\;\;\;whenever\;\;(i,j)\in J
\end{equation}
for sufficiently small $|z|>0$, where ${\mathcal G}_{i,j}$ is a short-hand notation 
for ${\mathcal G}_{A_{i,j}}^{(i,j)}$.
\end{Lemma}
{\it Proof.}
Using Proposition 6.1, we obtain
$$
A\omega_{i,j}(z)=\frac{1}{z}(\omega_{i,j}(z)-\zeta_{i,j})+\mathcal{R}_{A}(z)\omega_{i,j}(z)+
\alpha_{i,j}Df_{i,j}(\alpha_{i,j}^{2}z)\Omega_{i,j}
$$
for any $(j,j)\in J$. This equation is very similar to (4.4), except for the term involving the difference quotient.  
However, the projection $P_{i,j}$ maps the vector $\Omega_{i,j}$ to zero and thus
$$
A_{i,j}\omega_{i,j}(z)= \frac{1}{z}(\omega_{i,j}(z)-\zeta_{i,j})+\mathcal{R}_{A}(z)\omega_{i,j}(z)
$$
This leads to the equation
$$
z\omega_{i,j}(z)= \left(\frac{1}{z}+{\mathcal R}_{A}(z)-A_{i,j}\right)^{-1}\zeta_{i,j}
$$
for small $|z|$, as in the proof of Lemma 4.2, which proves (6.10) 
since $\langle \omega_{i,j}(z), \zeta_{i,j}\rangle = 1$. 
\hfill $\blacksquare$\\

It is not hard to see that a result similar to Lemma 5.1 can be established for conjugate states $\psi_{i,j}$, which
is stated below. Note that we can use the same $C(z)$ and its multiplicative inverse $B(z)$ in computations 
involving moments of $A$ in the state $\varphi$ as well as in its 
conjugate states $\psi_{i,j}$.
\begin{Lemma}
An ${\mathcal I}$-valued power series $\;\mathcal{R}_{A}(z)=\sum_{n=1}^{\infty}c_{n}z^{n-1}$ 
converging in the norm topology in a neighborhood of zero satisfies the equation (6.10) if and only if
\begin{equation}\tag{6.11}
\sum_{k=1}^{m}\sum_{n_1+\ldots +n_k=m-k}
\varphi_{i,j}(b_{n_1}A_{i,j}b_{n_2}\ldots b_{n_{k-1}}A_{i,j}b_{n_k})=0
\end{equation}
for all $m\geq 2$, where we assume that $n_1, \ldots , n_k$ are non-negative integers and 
where  $B(z)=\sum_{n=0}^{\infty}b_{n}z^{n+1}$ is the multiplicative inverse of $C(z)=1/z+\mathcal{R}(z)$.
\end{Lemma}
{\it Proof.}
The proof is similar to that of Lemma 5.1.
\hfill $\blacksquare$\\

\begin{Theorem}
Under the assumptions of Theorem 5.1, the series of the form
\begin{equation}\tag{6.12}
\mathcal{R}_{A}(z)=\sum_{i,j}\mathcal{R}_{i,j}(z)
\end{equation}
is an operatorial R-transform of the $\psi_{i,j}$-distribution of $A_{i,j}$ for any $(i,j)\in J$.
\end{Theorem}
{\it Proof.}
Lemma 6.2 enables us to carry out the same approximation procedure for arbitrary arrays of strongly matricially 
free random variables in terms of operators of Toeplitz type as in the proof of Theorem 5.1, 
which leads us to the desired conclusion. 
\hfill $\blacksquare$\\

It turns out that if we require an operatorial R-transform to satisfy (5.6) for the state $\varphi$ and 
(6.11) for all states $\psi_{i,j}$, then it is unique, and therefore it must be the matricial R-transform. 
In particular, in the Fock-space framework, these are the states defined by $\Omega$ and the array $(\Omega_{i,j})$.
Note that there is no need to study conjugate states of orders higher than two. 
On the Fock-space level (or, on the level of Hilbert space representations), 
the reason is that the action of strongly matricially free operators depends only on the type of the
first vector in any simple tensor which defines a state.

\begin{Theorem}
Let ${\mathcal T}_{i,j}(z)=T_{i,j}(z)1_{i,j}$ for any $(i,j)\in J$, where 
each $T_{i,j}(z)$ is a power series converging in some neighborhood of zero. Then there 
exists a unique series
\begin{equation}\tag{6.13}
\mathcal{R}_{A}(z)=\sum_{i,j}{\mathcal T}_{i,j}(z)
\end{equation}
which is an operatorial R-transform of the distribution of $A$ in the state $\varphi$ 
and of the distribution of $A_{i,j}$ in the conjugate state $\psi_{i,j}$ for any $(i,j)\in J$.
\end{Theorem}
{\it Proof.}
Existence of an operatorial R-transform of this form is established by Theorem 6.1. 
Conditions (5.6) and (6.11) uniquely determine $\varphi(b_{m})$ and $\varphi_{j}(b_{m})$ 
for all $m$'s and $j\in \{1,2\}$ since they are simple recursions. 
In the case of $\varphi$, we showed in Remark 5.3 how to solve it,
but the recursions for $\varphi_1$ and $\varphi_2$ are treated in a similar manner.
Now, using
$$
\varphi(\sum_{i,j}\beta_{i,j}q_{i,j})=\beta_{1,1},\;\;
\varphi_{1}(\sum_{i,j}\beta_{i,j}q_{i,j})=\beta_{2,1},\;\;
\varphi_{2}(\sum_{i,j}\beta_{i,j}q_{i,j})=\beta_{1,2},
$$
where pojections $(q_{i,j})$ are given by Remark 2.1,
we uniquely determine 
${\mathcal T}_{1,1}+{\mathcal T}_{2,2}, {\mathcal T}_{1,2}+{\mathcal T}_{2,2}, {\mathcal T}_{2,1}+{\mathcal T}_{1,1}$,
and thus also ${\mathcal T}_{2,1}+{\mathcal T}_{1,2}$, which implies that 
$B_{A}(z)$ is uniquely determined, and consequently, its multiplicative inverse $C_{A}(z)$ 
and thus ${\mathcal R}_{A}(z)$ in the orthogonal form
$$
\mathcal{R}_{A}(z)=\sum_{i,j}{\mathcal Q}_{i,j}(z)
$$
where ${\mathcal Q}_{i,j}(z)=Q_{i,j}(z)q_{i,j}$ for any $(i,j)\in J$, with $Q_{i,j}$'s
being ${\mathcal I}$-valued power series derived in the proof of Proposition 5.2, 
and where we identify $q_{i,j}$'s with elements of ${\mathcal I}$ defined as in Remark 2.1.
This completes the proof of uniqueness.
\hfill$\blacksquare$\\

\section{Cumulants}

In this section we return to the study of the moments of the strongly matricially free convolution 
in the state $\varphi$, which can be identified with the collection of moments 
\begin{equation}\tag{7.1}
\{\varphi(A^{m}): m\geq 0\}.
\end{equation}
Our results on the matricial R-transform indicate that there should 
exist a natural combinatorial approach which would allow us to express these moments in terms of free cumulants rather than to introduce new cumulants. Having in mind Theorem 3.2, this would
also give an alternative approach to the c-free convolution [4] and computations of its moments without
using c-free cumulants.

If $\mu$ is the distribution of a random variable with moments $\{M_{\mu}(m):m \geq 0\}$, then
\begin{equation}\tag{7.2}
M_{\mu}(m)=\sum_{\pi=\{\pi_1, \ldots , \pi_r\}\in \mathcal{NC}_{m}}r_{\mu}(|\pi_1|)\ldots r_{\mu}(|\pi_r|)
\end{equation}
where the numbers $(r_{\mu}(n))_{n\geq 1}$ are the free cumulants of $\mu$ [19]. 

Let us apply free cumulants to the combinatorics of the strongly matricially free convolution
of $(\mu_{i,j})$. Let $F_{n}(\pi)$ be the set of all mappings $f\in {\mathcal B}(\pi)\rightarrow \{1, 2, \ldots , n\}$
called {\it colorings} of the blocks ${\mathcal B}(\pi)$ of $\pi$, where $n\in {\mathbb N}$.
Then the pair $(\pi,f)$ plays the role of a {\it colored partition} with blocks denoted
\begin{equation}\tag{7.3}
{\mathcal B}(\pi,f)=\{(\pi_1,f),(\pi_2,f), \ldots , (\pi_r,f)\}
\end{equation}
to which we shall assign free cumulants in a suitable `matricial' way.

In the case of the free convolution, 
if we are given distributions $\mu_{1}, \mu_2$ whose
free cumulants of order $k$ are given by $(r_{j}(k))_{1\leq j \leq 2}$, then we can express the moments of
$\mu=\mu_1\boxplus \mu_2$ as 
\begin{equation}\tag{7.4}
M_{\mu}(m)=
\sum_{(\pi,f)=\{(\pi_1,f), \ldots , (\pi_r,f)\}\in \mathcal{NC}_{m}(2)}
r(\pi_1,f)\ldots r(\pi_r,f)
\end{equation}
where
\begin{equation}\tag{7.5}
r(\pi_k,f)=r_{j}(|\pi_k|)\;\;{\rm whenever}\;\;f(\pi_k)=j,
\end{equation}
and where $\mathcal{NC}_{m}(2)$ denotes the set of non-crossing partitions of the set 
$\{1,2, \ldots , m\}$ which are colored by the set $\{1,2\}$. 

We would like to use a similar formalism to describe the strongly matricially free convolution.
In contrast to (7.5), where cumulants are assigned to blocks independently 
of other blocks, we shall now make them depend on the colors 
of their outer blocks. For that purpose, instead of a tuple of cumulants of order $k$, 
we shall consider an array of free cumulants $(r_{i,j}(k))$ for each $k$, where $r_{i,j}(k)$ is the free 
cumulant of order $k$ of the distribution $\mu_{i,j}$, where $(i,j)\in J$. Clearly, this
is in agreement with the matricial formalism on which our approach is based. 

We assign cumulants to colored blocks of $(\pi,f)\in \mathcal{NC}_{m}(2)$ as follows.
If the block $(\pi_k,f)$ and all its outer blocks are colored by $j$, we assign to it the diagonal 
cumulant 
\begin{equation}\tag{7.6}
r(\pi_k,f)=r_{j,j}(|\pi_k|).
\end{equation}
On the other hand, if the block $(\pi_k,f)$ is colored by $i$ and it has 
an outer block colored by another color, then we assign to $(\pi_k,f)$ the off-diagonal cumulant
\begin{equation}\tag{7.7}
r(\pi_k,f)=r_{i,j}(|\pi_k|)
\end{equation}
where $j$ is the color of the deepest outer block of $\pi_k$ which has that property ($j\neq i$).
In the case when some $\mu_{i,j}$ does not appear in the considered array $(\mu_{i,j})$, 
we set $r_{i,j}(k)=0$ for any $k$ (one can also formally set $\mu_{i,j}=\delta_{0}$, which makes all cumulants vanish).
Note that we do not require this outer block to be the nearest outer block of $\pi_k$ as in [12-13], 
where the combinatorics is slightly simpler since it is based on Gaussian operators, but we 
do require it to be its nearest outer block which is {\it differently colored}. 
The examples of Fig.1 (especially the last one) show this feature in more detail. 

Therefore, if $\pi_k$ is a block colored by $j$ and all its outer blocks are colored by $j$, 
we assign to it the pair $(j,j)$. If $\pi_k$ is colored by $i$ and it 
has outer blocks wich are differently colored, then we assign to it
the pair $(i,j)$, where $j$ is the color of the deepest blocks among those.
We will then say that $\pi_k$ is {\it labelled} by $(j,j)$ or $(i,j)$, respectively.
The labelling $(i,j)$ induced by the coloring will be called {\it $J$-admissible} if $(i,j)\in J$. 

In order to eliminate labellings which are not $J$-admissible from our formulas, we will denote by 
$\mathcal{NC}_{m}(2,J)$ the subset of $\mathcal{NC}_{m}(2)$ consisting of those 
partitions which induce $J$-admissible labellings. For instance, if $(1,2)\notin J$, then 
all colored partitions in Fig.1 but the first one are not $J$-admissible.

\begin{Definition}
{\rm
Let $\mu$ be the distribution of $A$ given by (3.1) and let $(\mu_{i,j})$ be the corresponding array of 
distributions. By the {\it partitioned colored cumulant} corresponding to the 
colored partition $(\pi,f)\in \mathcal{NC}_{m}(2)$ we understand the product
\begin{equation}\tag{7.8}
r_{\mu}[\pi,f]=r(\pi_1,f)\ldots r(\pi_r,f)
\end{equation}
where $\pi=\{\pi_1, \pi_2, \ldots, \pi_r\}$ and free cumulants are assigned to its colored blocks 
according to (7.6)-(7.7).} 
\end{Definition}

As we already remarked in Section 4, 
it suffices to consider Toeplitz operators to compute the moments of $\boxplus_{i,j}\mu_{i,j}$. 
We use these operators to derive the moment-cumulant formula given below.
In contrast to the free or the c-free case, this formula is not the definition of cumulants since 
it expresses the moments of the sum of random variables in terms of free cumulants 
associated with marginal laws.

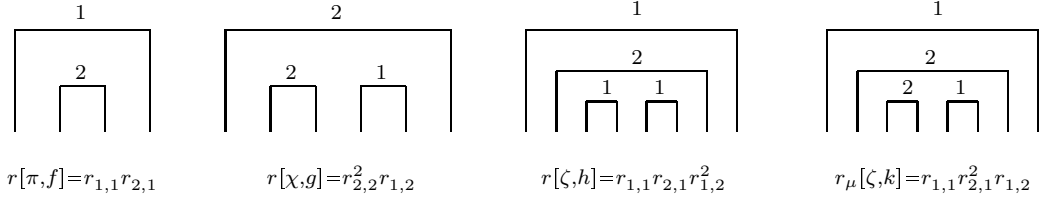
\begin{figure}
\unitlength=1mm
\special{em:linewidth 0.4pt}
\linethickness{0.4pt}
\begin{picture}(120.00,30.00)(24.00,2.00)

\put(17.00,10.00){\line(0,1){13.50}}
\put(23.00,10.00){\line(0,1){6.00}}
\put(29.00,10.00){\line(0,1){6.00}}
\put(35.00,10.00){\line(0,1){13.50}}

\put(25.00,25.00){$\scriptstyle{1}$}
\put(25.00,17.00){$\scriptstyle{2}$}

\put(17.00,23.50){\line(1,0){18.00}}
\put(23.00,16.00){\line(1,0){6.00}}
\put(16.00,3.00){$\scriptstyle{r[\pi,f]=r_{1,1}^{}r_{2,1}^{}}$}


\put(45.00,10.00){\line(0,1){13.50}}
\put(51.00,10.00){\line(0,1){6.00}}
\put(57.00,10.00){\line(0,1){6.00}}
\put(63.00,10.00){\line(0,1){6.00}}
\put(69.00,10.00){\line(0,1){6.00}}
\put(75.00,10.00){\line(0,1){13.50}}

\put(53.00,17.00){$\scriptstyle{2}$}
\put(65.00,17.00){$\scriptstyle{1}$}
\put(59.00,25.00){$\scriptstyle{2}$}

\put(45.00,23.50){\line(1,0){30.00}}
\put(51.00,16.00){\line(1,0){6.00}}
\put(63.00,16.00){\line(1,0){6.00}}
\put(50.50,3.00){$\scriptstyle{r[\chi,g]=r_{2,2}^{2}r_{1,2}^{}}$}


\put(85.00,10.00){\line(0,1){13.50}}

\put(89.00,10.00){\line(0,1){8.00}}
\put(93.00,10.00){\line(0,1){4.00}}
\put(97.00,10.00){\line(0,1){4.00}}
\put(101.00,10.00){\line(0,1){4.00}}
\put(105.00,10.00){\line(0,1){4.00}}
\put(109.00,10.00){\line(0,1){8.00}}

\put(113.00,10.00){\line(0,1){13.50}}

\put(95.00,15.00){$\scriptstyle{1}$}
\put(102.00,15.00){$\scriptstyle{1}$}
\put(99.00,19.00){$\scriptstyle{2}$}
\put(99.00,25.50){$\scriptstyle{1}$}

\put(85.00,23.50){\line(1,0){28.00}}
\put(89.00,18.00){\line(1,0){20.00}}
\put(93.00,14.00){\line(1,0){4.00}}
\put(101.00,14.00){\line(1,0){4.00}}
\put(87.00,3.00){$\scriptstyle{r[\zeta,h]=r_{1,1}^{}r_{2,1}^{}r_{1,2}^{2}}$}

\put(125.00,10.00){\line(0,1){13.50}}

\put(129.00,10.00){\line(0,1){8.00}}
\put(133.00,10.00){\line(0,1){4.00}}
\put(137.00,10.00){\line(0,1){4.00}}
\put(141.00,10.00){\line(0,1){4.00}}
\put(145.00,10.00){\line(0,1){4.00}}
\put(149.00,10.00){\line(0,1){8.00}}

\put(153.00,10.00){\line(0,1){13.50}}

\put(135.00,15.00){$\scriptstyle{2}$}
\put(142.00,15.00){$\scriptstyle{1}$}
\put(138.00,19.00){$\scriptstyle{2}$}
\put(139.00,25.50){$\scriptstyle{1}$}

\put(125.00,23.50){\line(1,0){28.00}}
\put(129.00,18.00){\line(1,0){20.00}}
\put(133.00,14.00){\line(1,0){4.00}}
\put(141.00,14.00){\line(1,0){4.00}}
\put(126.00,3.00)
{$\scriptstyle{r_{\mu}[\zeta,k]=r_{1,1}^{}r_{2,1}^{2}r_{1,2}^{}}$}

\end{picture}
\caption{Partitioned colored cumulants}
\end{figure}

\begin{Lemma}
If $\mu$ is the $\varphi$--distribution of the sum $A=\sum_{i,j}a_{i,j}$ of strongly matricially free random variables
and $\mu_{i,j}$ is the $\varphi_{i,j}$--distribution of $a_{i,j}$, then 
\begin{equation}\tag{7.9}
M_{\mu}(m)=\sum_{(\pi, f)\in \mathcal{NC}_{m}(2,J)} r_{\mu}[\pi,f]
\end{equation}
where the number $r_{\mu}[\pi,f]$ is the partitioned colored cumulant corresponding to $(\pi, f)$.
\end{Lemma}
{\it Proof.}
Without loss of generality, we can choose $(a_{i,j})$ to be a two-dimensional array of Toeplitz operators of the form
(4.1), where 
\begin{equation*}
f_{i,j}(\ell_{i,j}^{*})=\sum_{k=1}^{m}s_{i,j}(k)(\ell_{i,j}^{*})^{k-1},
\end{equation*}
for some positive integer $m$, where we adopt the convention that $(\ell_{i,j}^{*})^{0}=1_{i,j}$
for any $(i,j)\in J$. Clearly, we obtain
$$
\varphi(A^{m})=\sum_{((i_1,j_1), \ldots,(i_m,j_m))\in \Delta_{m}}
\varphi(a_{i_1,j_1}\ldots a_{i_m,j_m})
$$
where $\Delta_{m}$ is the subset of $J^{m}$, for which there exist operators 
\begin{equation*}
b_k\in \{\ell_{i_k,j_k}, (\ell_{i_k,j_k}^{*})^{p-1}; p\geq 1\}, \;\;1\leq k \leq m,
\end{equation*}
such that
\begin{equation}\tag{7.10}
\varphi(b_1b_2\ldots b_m)\neq 0.
\end{equation}
We claim that the set of tuples $(b_1,b_2, \ldots , b_m)$ which satisfy these conditions is 
in 1-1 correspondence with $\mathcal{NC}_{m}(2,J)$.

Namely, to each $(\pi,f)\in \mathcal{NC}_{m}(2,J)$ 
we assign a unique tuple $(b_1,b_2, \ldots , b_m)$ according to the following rules:
\begin{enumerate}
\item[(a)]
if $k$ corresponds to the first leg of a $p$-block labelled by $(i,j)$,
we assign to it $b_k=(\ell_{i,j}^{*})^{p-1}$,
\item[(b)]
if $k$ corresponds to any but the first leg of a $p$-block labelled by $(i,j)$, 
we assign to it $b_k=\ell_{i,j}$,
\end{enumerate}
where by a $p$-block we understand a block consisting of $p$ elements. We claim that in this case
(7.10) holds. In fact, if $\{k, \ldots, k+p-1\}$ is a block for some $k$ and $p$, 
we have
\begin{equation}\tag{7.11}
(b_{k},b_{k+1}, \ldots , b_{k+p-1})=((\ell_{i,j}^{*})^{p-1}, \ell_{i,j}, \ldots , \ell_{i,j})
\end{equation} 
where 
\begin{equation}\tag{7.12}
(i,j)=(i_k,j_k)=\ldots=(i_{k+p-1},j_{k+p-1}).
\end{equation}
Using (2.8), we can see that
the associated product produces a power of $\alpha_{i,j}$, which reduces the computation of the given moment 
to a moment of order $s<m$ and we can apply a similar procedure to this reduced moment. 
The whole procedure finally gives a product of $\alpha_{i,j}$'s and thus (7.10) holds 
since $\alpha_{i,j}\neq 0$ for any $(i,j)\in J$ by assumption.

Conversely, suppose that a product $b_1b_2\ldots b_m$ of the considered type, with a nonvanishing moment in the state $\varphi$, is given.
We would like to find a unique $(\pi,f)\in \mathcal{NC}_{m}(2,J)$ 
associated to this product according to the above rules.
This will be done by induction with respect to $m$. First, observe that if $m=1$, then 
the non-vanishing moment must be of type $\varphi(1_{j,j})$, where $(j,j)\in J$ 
by Proposition 2.1(2), which corresponds to the partition consisting of a single 
1-block colored by $j$ with a $J$-admissible labelling $(j,j)$.
If $m>1$ and $b_k=1_{i_k,j_k}$ for all $1\leq k \leq m$, then we must have $i_k=j_k$ and $(j_k,j_k)\in J$ 
for all $k$, which forces $(\pi,f)$ to consist of $1$-blocks only (they all have $J$-admissible labellings) 
and thus they must be colored by $i_1,i_2, \ldots, i_m$, respectively.

If $m>1$ and not all operators are units, there exists $1\leq k \leq m$ such that (7.11) and (7.12) hold
for some $p\geq 1$ since otherwise either all $b_{i}$'s would be creation operators, 
or all annihilation operators would be followed by other annihilation operators or 'wrong' creation operators, 
which would make the moment vanish. From among tuples of operators of this type, it is convenient to choose the one with largest $k$, 
which implies that $b_{k+p-1}=b_m$ or $b_{k+p}$ is a creation operator with $i_{k+p}=j$.
Again, using the relation (2.8), we can see that such a product contributes a power of $\alpha_{i,j}$.
This reduces the moment to a non-vanishing moment of order $s<m$, to which there corresponds a unique $\pi'\in \mathcal{NC}_{s}(2,J)$
with a $J$-admissible labelling by the inductive assumption. 
Now, there exists a unique $\pi\in \mathcal{NC}_{m}(2,J)$ obtained from $\pi'$ by inserting the block $\{k,\ldots , k+p-1\}$ colored by 
$i$ in between $k-1$ and $k+p$. Of course, if $k+p-1=m$, this block is separated from $\pi'$ and we must have $i=j$ 
since off-diagonal creation operators kill the vacuum, which uniquely determines the $J$-admissible labelling of the considered block.
Otherwise, the inserted block is inner with respect to a block containing $k+p$ and its $J$-admissible 
labelling $(i,j)$ is uniquely determined by the color of that block.

Recall that the contribution from the product of operators given in (7.11), to which we 
associate a $p$-block, should be multiplied by the coefficient $s_{i,j}(p)$, 
which gives the contribution
$$
s_{i,j}(p)\alpha_{i,j}^{2(p-1)}=r_{i,j}(p)
$$
from this block to the moment $\varphi(a_{i_1,j_1}a_{i_2,j_2}\ldots a_{i_m,j_m})$. 
The product of these contributions gives the product of the cumulants corresponding to all blocks of $(\pi,f)$. 
Therefore, the contribution from the product of operators assigned to any $(\pi,f)\in {\mathcal NC}_{m}(2,J)$ is 
$$
r(\pi_1,f)\ldots r(\pi_r,f),
$$
which completes the proof of the combinatorial formula for the moments of $A$.
\hfill $\blacksquare$

\begin{Example}
{\rm To colored non-crossing pair-partitions shown in Fig.1 we assign partitioned colored cumulants.
For instance, we have $r[\pi,f]=r(\pi_1,f)r(\pi_2,f)$, where
$\pi_1=\{1,4\}$ and $\pi_2=\{2,3\}$. Since all free cumulants involved are of order 2, we write $r_{i,j}$ instead of $r_{i,j}(2)$, which simplifies the notation. If we collect all colorings for the considered partitions, we 
obtain the sums of partitioned colored cumulants over all colorings:
\begin{eqnarray*}
r[\pi]&=& r_{1,1}(r_{1,1}+r_{2,1}) + r_{2,2}(r_{2,2}+r_{1,2}),\\
r[\chi]&=& r_{1,1}(r_{1,1}+r_{2,1})^{2}+r_{2,2}(r_{2,2}+r_{1,2})^{2}\\
r[\zeta]&=&r^{2}_{1,1}(r_{1,1}^{}+r_{2,1}^{})^{2}+r_{1,1}^{}r_{2,1}^{}r_{1,2}^{2}\\
&& +
r_{2,2}^{2}(r_{2,2}^{}+r_{1,2}^{})^{2} + r_{2,2}^{}r_{1,2}^{}r_{2,1}^{2}
\end{eqnarray*}
Since $\zeta$ is a partition of depth three, we can observe the influence of  
the definition of {\it strongly} matricially free r.v. 
Namely, if blocks $\zeta_1=\{1,8\}$ and $\zeta_2=\{2,7\}$ 
are differently colored, then blocks which are inner with respect to $\zeta_2$, namely $\zeta_3=\{3,4\}$ and
$\zeta_4=\{5,6\}$ have to be colored differently than $\zeta_2$ due to the fact that 
we consider the case of strongly matricially free random variables.
If we set $r_{1,2}=r_{1,1}=r_1$ and $r_{2,1}=r_{2,2}=r_2$ 
(row-identically distributed square array), we obtain 
\begin{eqnarray*}
r[\pi]&=& r_{1}^{2}+2r_{1}r_{2}+r_{2}^{2},\\
r[\chi]&=& r_{1}^{3}+3r_{1}^2r_{2}+3r_{1}r_{2}^2+r_{2}^{3},\\
r[\zeta]&=& r_{1}^{4}+3r_{1}^3r_{2}+2r_{1}^{2}r_{2}^{2}+3r_{1}r_{2}^3+r_{2}^{4},
\end{eqnarray*}
which is the contribution from partitions $\pi$, $\chi$ and $\zeta$
to the moments of $\mu_1\boxplus \mu_2$. In turn, if we set $r_{2,1}=r_{2,2}=r_2$, $r_{1,1}=r_1$ and $r_{1,2}=0$ (row-identically distributed lower-triangular 
array), we obtain
\begin{eqnarray*}
r[\pi]&=& r_{1}^{2}+r_{1}r_{2}+r_{2}^{2},\\
r[\chi]&=& r_{1}^{3}+r_{1}^2r_{2}+r_{1}r_{2}^2+r_{2}^{3},\\
r[\zeta]&=& r_{1}^{4}+2r_{1}^3r_{2}+r_{1}^{2}r_{2}^{2}+r_{2}^{4},
\end{eqnarray*}
which gives the contribution from $\pi$, $\chi$ and $\zeta$ to 
the moments of $\mu_1\vartriangleright \mu_2$.}
\end{Example}

\begin{Example}
{\rm
The lowest order moments of $A$ are expressed in terms of free cumulants 
of the distributions $\mu_{i,j}$ as follows:
\begin{eqnarray*}
M_{\mu}(1)&=& r_{1,1}(1)+r_{2,2}(1),\\
M_{\mu}(2)&=& r_{1,1}(2)+r_{2,2}(2) + 
(r_{1,1}(1)+r_{2,2}(1))^{2},\\
M_{\mu}(3)&=&r_{1,1}(3)+r_{2,2}(3) + 2(r_{1,1}(2)+r_{2,2}(2))(r_{1,1}(1)+r_{2,2}(1))\\
&+&r_{1,1}(2)(r_{1,1}(1)+r_{2,1}(1))+r_{2,2}(2)(r_{2,2}(1)+r_{1,2}(1))\\
&+&(r_{1,1}(1)+r_{2,2}(1))^{3}.
\end{eqnarray*}
As in the case of moments of arbitrary orders, these moments agree with the moments of $\mu_1\boxplus \mu_2$ 
if the array $(a_{i,j})$ is square and row-identically distributed. 
In this case the moments of $A$ can be expressed in terms of free cumulants of $\mu_1\boxplus \mu_2$.
In turn, if that array is lower-triangular and row-identically distributed, these moments
agree with those of $\mu_1\vartriangleright \mu_2$. However, as already $M_{\mu}(3)$ demonstrates, 
the partitioned colored cumulants which appear in our moment-cumulant formula
cannot be expressed in terms of free cumulants of free convolutions of distributions.
Roughly speaking, this effect shows how much the shape of the lower-triangular array
affects the additivity of the considered transforms.}
\end{Example}

\section{Cauchy transforms}
Continuing the considerations of Section 7, we shall study the Cauchy transforms of the 
`commutative' $\varphi$-distribution of $A$ and express them in terms of their R-transforms. 
One should remember that these are Cauchy transforms with complex arguments rather than
those with operatorial arguments which provided the right framework 
for the linearization property of the R-transform. 

As a consequence of Proposition 4.1 and Lemma 7.1, the $\varphi$-distribution of $A$ is uniquely determined by
the array of free cumulants of distributions $(\mu_{i,j})$. The way we defined partitioned colored cumulants indicates
that the moments of $A$ under $\varphi$ exhibit the property of `diagonal subordination', 
which is reflected by the fact that all covering blocks of non-crossing partitions which appear
in the formula for the moments of $A$ in terms of free cumulants 
are labelled by diagonal pairs $(j,j)$, where $j\in \{1,2\}$, and all inner blocks have 
colorings which are `matricially subordinate' to their closest outer blocks.

More precisely, the first level (counting from the top) of each non-crossing 
partition (i.e. all covering blocks, including singletons) is labelled by some $(j,j)$, but
all blocks which are inner with respect to a given block $\pi_k$ labelled by $(j,j)$
and which form nearest inner-outer pairs with that block 
(thus, all blocks at the second level counting from the top) 
are `subordinate' to $\pi_k$ and the corresponding coloring in the sense that
its labelling must be $(i,j)$ for some $i$.
A similar `subordination' holds on the remaining levels.

\begin{figure}
\unitlength=1mm
\special{em.linewidth 0.5pt}
\linethickness{0.5pt}
\begin{picture}(140.00,60.00)(-25.00,5.00)

\put(10.00,10.00){\line(1,2){5.00}}
\put(20.00,10.00){\line(-1,2){5.00}}
\put(30.00,10.00){\line(1,2){5.00}}
\put(40.00,10.00){\line(-1,2){5.00}}

\put(50.00,10.00){\line(1,2){5.00}}
\put(60.00,10.00){\line(-1,2){5.00}}
\put(70.00,10.00){\line(1,2){5.00}}
\put(80.00,10.00){\line(-1,2){5.00}}

\put(15.00,20.00){\line(2,3){10.00}}
\put(35.00,20.00){\line(-2,3){10.00}}
\put(55.00,20.00){\line(2,3){10.00}}
\put(75.00,20.00){\line(-2,3){10.00}}

\put(25.00,35.00){\line(4,3){20.00}}
\put(65.00,35.00){\line(-4,3){20.00}}

\put(10.00,10.00){\circle*{1.50}}
\put(20.00,10.00){\circle*{1.50}}
\put(30.00,10.00){\circle*{1.50}}
\put(40.00,10.00){\circle*{1.50}}
\put(50.00,10.00){\circle*{1.50}}
\put(60.00,10.00){\circle*{1.50}}
\put(70.00,10.00){\circle*{1.50}}
\put(80.00,10.00){\circle*{1.50}}

\put(10.00,10.00){\line(-1,-2){2.00}}
\put(10.00,10.00){\line(1,-2){2.00}}
\put(20.00,10.00){\line(-1,-2){2.00}}
\put(20.00,10.00){\line(1,-2){2.00}}
\put(30.00,10.00){\line(-1,-2){2.00}}
\put(30.00,10.00){\line(1,-2){2.00}}
\put(40.00,10.00){\line(-1,-2){2.00}}
\put(40.00,10.00){\line(1,-2){2.00}}

\put(50.00,10.00){\line(-1,-2){2.00}}
\put(50.00,10.00){\line(1,-2){2.00}}
\put(60.00,10.00){\line(-1,-2){2.00}}
\put(60.00,10.00){\line(1,-2){2.00}}
\put(70.00,10.00){\line(-1,-2){2.00}}
\put(70.00,10.00){\line(1,-2){2.00}}
\put(80.00,10.00){\line(-1,-2){2.00}}
\put(80.00,10.00){\line(1,-2){2.00}}

\put(15.00,20.00){\circle*{1.50}}
\put(35.00,20.00){\circle*{1.50}}
\put(55.00,20.00){\circle*{1.50}}
\put(75.00,20.00){\circle*{1.50}}

\put(25.00,35.00){\circle*{1.50}}
\put(65.00,35.00){\circle*{1.50}}

\put(45.00,50.00){\circle*{1.50}}

\put(18.00,35.00){\footnotesize $\mu_{1,1}$}
\put(67.00,35.00){\footnotesize $\mu_{2,2}$}

\put(8.00,20.00){\footnotesize $\mu_{1,1}$}
\put(27.00,20.00){\footnotesize $\mu_{2,1}$}
\put(58.00,20.00){\footnotesize $\mu_{1,2}$}
\put(77.00,20.00){\footnotesize $\mu_{2,2}$}

\put(3.00,10.00){\footnotesize $\mu_{1,1}$}
\put(13.00,10.00){\footnotesize $\mu_{2,1}$}
\put(22.50,10.00){\footnotesize $\mu_{1,2}$}
\put(33.00,10.00){\footnotesize $\mu_{2,1}$}
\put(51.00,10.00){\footnotesize $\mu_{1,2}$}
\put(61.00,10.00){\footnotesize $\mu_{2,1}$}
\put(72.00,10.00){\footnotesize $\mu_{1,2}$}
\put(81.00,10.00){\footnotesize $\mu_{2,2}$}

\end{picture}
\caption{Binary subordination of measures}
\end{figure}
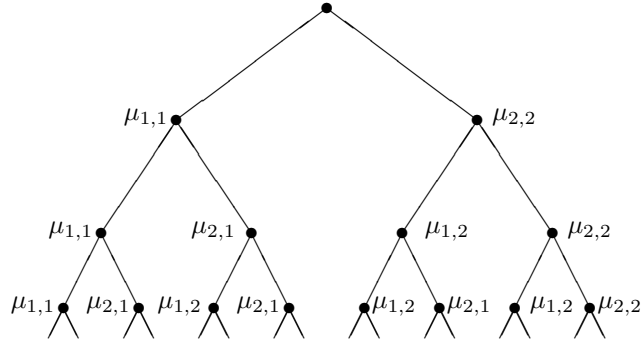

Therefore, it is natural to introduce a mapping which assigns 
to each distribution $\mu_{i,j}$ the variable $A_{i,j}$ associated 
with the array $(\mu_{i,j}^{*})$ which is `subordinate' to
$\mu_{i,j}$ (to the top covering block we assign $r_{i,j}$).

\begin{Definition}
{\rm Let $(\mu_{i,j})$ be an array of distributions. By the {\it subordinate} arrays 
of distributions we understand arrays of the form
\begin{equation}\tag{8.1}
\mu_{1,1}^{*}=\left(
\begin{array}{cc}
\mu_{1,1} & \mu_{1,2}\\
\mu_{2,1} & \mu_{2,1}
\end{array}
\right),\;\;\;
\mu_{1,2}^{*}=\left(
\begin{array}{cc}
\mu_{1,2} & \mu_{1,2}\\
\mu_{2,1} & \mu_{2,1}
\end{array}
\right),
\end{equation}
\begin{equation}\tag{8.2}
\mu_{2,1}^{*}=\left(
\begin{array}{cc}
\mu_{1,2} & \mu_{1,2}\\
\mu_{2,1} & \mu_{2,1}
\end{array}
\right),\;\;\;
\mu_{2,2}^{*}=\left(
\begin{array}{cc}
\mu_{1,2} & \mu_{1,2}\\
\mu_{2,1} & \mu_{2,2}
\end{array}
\right).
\end{equation}
The corresponding Cauchy transform will be called {\it subordinate Cauchy transforms} and will be denoted
by $G_{i,j}^{*}$, respectively.}
\end{Definition}

Figure 2 shows how to interpret the subordination of distributions exhibited by these 
equations. It has the form of a binary tree, where the part of the tree that lies 
below any $\mu_{i,j}$ shows which measures have to be replaced in the original 
array in order to produce a formula for the corresponding moments of 
strongly matricially free random variables. 

More explicitly, these are Cauchy transforms of the sums
of strongly matricially free random variables with distributions given by $(\mu_{i,j}^{*})$.
It is easy to see that they satisfy the recursion given by equations
\begin{eqnarray*}
G_{j,j}^{*}(z)&=&
\frac{1}
{
z-
K_{j,j}(z)-K_{\bar{j},j}(z)
},
\;\;\; {\rm where}\; j\in \{1,2\}\\
G_{j,\bar{j}}^{*}(z)&=&
\frac{1}
{
z-
K_{j,\bar{j}}(z)
-
K_{\bar{j},j}(z)
}, \;\;\;{\rm where}\;
j\in \{1,2\},
\end{eqnarray*}
with $\overline{1}=2$ and $\overline{2}=1$, 
where $K_{i,j}(z)=R_{i,j}(G_{i,j}^{*}(z))$ for any $i,j\in \{1,2\}$.

\begin{Theorem}
Under the assumptions of Lemma 7.1, the Cauchy transform of the distribution of $A$ 
in the state $\varphi$ is given by
\begin{equation}\tag{8.3}
G_{A}(z)=\frac{1}
{
z-
\sum_{j}R_{j,j}^{}(G_{j,j}^{*}(z))
}
\end{equation}
where $(G_{i,j}^{*})$ is the associated array of subordinate Cauchy 
transforms.
\end{Theorem}
{\it Proof.}
The general combinatorics of the proof reminds that for conditional freeness. The most 
important difference is that we use free cumulants instead of c-free cumulants.
Let 
$$
M_{A}(z)=z^{-1}G_{A}(z^{-1})
$$
be the moment generating function associated with the $\varphi$-distribution of $A$.
Then, by Lemma 7.1, we have
$$
M_{A}(z)-1=\sum_{n=1}^{\infty}z^{n}\sum_{(\pi,f)\in \mathcal{NC}_{n}^{c}}r[\pi,f].
$$
On the other hand, if our assertion is supposed to hold, we must have
$$
M_{A}(z)-1=\sum_{j=1}^{2}K_{j,j}(z^{-1})zM_{A}(z)
$$
where 
\begin{eqnarray*}
K_{j,j}(z^{-1})&=&\sum_{m=1}^{\infty}r_{j,j}(m)\sum_{n_1,\ldots ,n_{m-1}=0}^{\infty}
(zM_{j,j}^{*}(z))^{n_1}\ldots (zM_{j,j}^{*}(z))^{n_{m-1}}
\end{eqnarray*}
which leads to the equation
$$
M_{A}(z)-1=\sum_{j=1}^{2}\sum_{m=1}^{\infty}r_{j,j}(m)\sum_{n_1,\ldots ,n_{m}=0}^{\infty}
m_{j,j}^{*}(n_1)\ldots m_{j,j}^{*}(n_{m-1})m_{A}(n_m)z^{n_1+\ldots +n_{m}+m}.
$$
Therefore, the coefficient standing by $z^n$ in this expression is 
given by
$$
c_{A}(n)=\sum_{j=1}^{2}\sum_{m=1}^{n}\sum_{n_1+\ldots +n_{m}=n-m}
r_{j,j}(m)m_{j,j}^{*}(n_1)\ldots m_{j,j}^{*}(n_{m-1})m_{A}(n_m)
$$
We need to compare this number with $m_{A}(n)$ expressed in terms of free cumulants:
$$
m_{A}(n)=\sum_{(\pi,f)\in \mathcal{NC}_{n}^{c}}r[\pi,f]
$$
If we express each moment $m_{j,j}^{*}(n_k)$ in the above formula for $c_{A}(n)$ in terms of free cumulants, 
we can argue that we obtain the same expression as that for $m_{A}(n)$. In fact, it suffices to interpret the product
$$
r_{j,j}(m)m_{j,j}^{*}(n_1)\ldots m_{j,j}^{*}(n_{m-1})m_{A}(n_m)
$$
as follows. The cumulant $r_{j,j}(m)$ is assigned to this block of some $(\pi,f)\in \mathcal{NC}_{n}^{c}$ 
which contains $1$. This block is assumed to have $m$ legs and to all of these legs but the first one we assign `subordinate' moments $m_{j,j}^{*}(n_1), \ldots , m_{j,j}^{*}(n_{m-1})$ (each of these moments
is a sum of paritioned colored cumulants associated with the partitions of the subinterval
lying between the given leg and its left neighbor). Here, numbers $n_1, \ldots , n_{m-1}$ are cardinalities of 
these subintervals. This part of the product corresponds to the first covering block and its inner blocks.
It remains to assign $m_{A}(n_m)$ to the remaining covering blocks and its inner blocks which has the same 
general form as $m_{A}(n)$ except that the cardinality of the considered interval is smaller. 
It can be seen that each partitioned cumulant is obtained in this fashion and it is obtained precisely 
once since the numbers $j$, $m$, cardinalities $n_1, \ldots , n_{m-1}$ of the above-mentioned subintervals, together
with $m-1$ colored non-crossing partitions of these subintervals, and a fixed product of free cumulants 
taken from the cumulant expression for $m_{A}(n_m)$ determine a unique partitioned cumulant of $(\pi,f)$. 
This completes the proof. \hfill $\blacksquare$\\

\begin{Corollary}
Under the assumptions of Theorem 8.1, the Cauchy transforms of the following 
convolutions are expressed in terms of R-transforms by formulas given below:
\begin{enumerate}
\item free convolution
\begin{equation*}
G_{\mu_1\boxplus\, \mu_2}(z)=\frac{1}{z-R_{\mu_1}(G_{\mu_1\boxplus \,\mu_2}(z))-R_{\mu_2}(G_{\mu_1\boxplus\, \mu_2}(z))}
\end{equation*}
\item monotone convolution
\begin{equation*}
G_{\mu_1\vartriangleright \mu_2}(z)=\frac{1}{z-R_{\mu_1}(G_{\mu_1\vartriangleright \mu_2}(z))-R_{\mu_2}(G_{\mu_2}(z))}
\end{equation*}
\item boolean convolution
\begin{equation*}
G_{\mu_1\uplus \mu_2}(z)=\frac{1}{z-R_{\mu_1}(G_{\mu_1}(z))-R_{\mu_2}(G_{\mu_2}(z))}
\end{equation*}
\item s-free convolution
\begin{equation*}
G_{\mu_1\boxright \,\mu_2}(z)=\frac{1}{z-R_{\mu_1}(G_{\mu_1\boxplus \mu_2}(z))}
\end{equation*}
\item orthogonal convolution
\begin{equation*}
G_{\mu_1\vdash \mu_1}(z)=\frac{1}{z-R_{\mu_1}(G_{\mu_1\vartriangleright \mu_2}(z))}
\end{equation*}
\end{enumerate}
\end{Corollary}
{\it Proof.} We shall present the proof of the first two formulas.
In the case of free convolution, we have 
$$
\mu_{i,j}^{*}=\left(
\begin{array}{cc}
\mu_{1} & \mu_{1}\\
\mu_{2} & \mu_{2}
\end{array}
\right),
$$
since distributions are row-identical and therefore
the array $\mu_{i,j}^{*}$ agrees with the original array $(\mu_{i,j})$. Hence,
$G_{i,j}^{*}=G_{\mu}$ for any $i,j$, and thus (8.3) easily gives the formula for the free convolution. 
In turn, in the case of monotone convolution, the array is lower-triangular and $r_{1,2}(k)=0$ for any $k$.
This implies that all blocks lying under blocks colored by $2$ must also be colored by $2$. Therefore, 
$\mu_{2,2}^{*}$ reduces to the one-dimensional array with $\mu_{2,2}=\mu_2$, whereas 
$$
\mu_{1,1}^{*}=\left(\begin{array}{cc}
\mu_{1,1} & \\
\mu_{2,1} & \mu_{2,1}
\end{array}
\right)=\left(\begin{array}{cc}
\mu_{1} & \\
\mu_{2} & \mu_{2}
\end{array}
\right)
$$
since blocks lying under the covering blocks colored by $1$ can be colored by $1$ or $2$. However, 
if such a block is colored by $2$, then the deepest block which is differently colored and is 
outer with respect to that block, must be colored by $1$. This gives $G_{1,1}^{*}=G_{\mu_1\vartriangleright \mu_2}$
and $G_{2,2}^{*}=G_{\mu_2}$, which gives the formula for the monotone convolution. The proofs of the 
remaining formulas are similar and are left to the reader.
\hfill $\blacksquare$

\begin{Remark}
{\rm The equations of Corollary 8.1 can also be proved without using Theorem 8.1. 
For instance, using formulas
$$
F_{\nu_1\vartriangleright\nu_2}=F_{\nu_1}\circ F_{\nu_2}\;\;{\rm and}\;\;F_{\nu}(z)=z-R_{\nu}(G_{\nu}(z)),
$$
where $F_{\nu}$ is the reciprocal Cauchy transform of $\nu$, we obtain
\begin{eqnarray*}
F_{\mu_1\vartriangleright \mu_2}(z)&=&F_{\mu_2}(z)-R_{\mu_1}(G_{\mu_1}(F_{\mu_2}(z))\\
&=& z-R_{\mu_2}(G_{\mu_2}(z))-R_{\mu_1}(G_{\mu_1\vartriangleright \mu_2}(z)),
\end{eqnarray*}
which is equivalent to equation (2). In a similar way, one can prove the remaining equations. 
}
\end{Remark}

Nevertheless, Theorem 8.1 gives a formula which relates Cauchy transforms of various 
convolutions to R-transforms of $\mu_1$ and $\mu_2$ in a unified manner. Let us also remark that
this formula can also be viewed as an approximation for Cauchy transforms of 
various convolutions.
In that respect, it is similar to formulas for the Cauchy transform of the 
free convolution of compactly supported measures given in [10].

\begin{Example}
{\rm Let us consider the example of a two-by-two square array in which the diagonal 
distributions are semicircle laws whereas the off-diagonal ones are point masses. 
Thus, let $R_{1,1}(z)=az$, $R_{2,2}(z)=dz$ and $R_{1,2}(z)=b$, $R_{2,1}(z)=c$. 
Then, using (8.3), we obtain
$$
G_{\mu}(z)=\frac{1}{z-aG_{1,1}^{*}(z)-dG_{2,2}^{*}(z)},
$$
where $\mu=\boxplus_{i,j}\mu_{i,j}$ and the diagonal subordinate Cauchy transforms correspond to shifted semicircle laws 
since they satisfy the equations
$$
G_{1,1}^{*}(z)=\frac{1}{z-c-aG_{1,1}^{*}(z)}\;\;\;
{\rm and}\;\;\;G_{2,2}^{*}(z)=\frac{1}{z-b-dG_{2,2}^{*}(z)}.
$$
In the general case, the analytic formula for $G_{\mu}$ is rather complicated, but if $a=d$ and $c=b$, then
$G_{1,1}^{*}=G_{2,2}^{*}$ and thus
$$
G_{\mu}(z)=\frac{1}{z-2aG_{1,1}^{*}(z)}=\frac{b-\sqrt{(b-z)^{2}-4a}}{4a+2bz-z^2},
$$
and therefore $\mu$ is the free Meixner distribution with the continuous density
$$
d\mu(x)=\frac{\sqrt{4a-(x-b)^{2}}}{\pi(4a+2bx-x^2)}
$$
on the interval $[b-2\sqrt{a}, b+2\sqrt{a}]$ and two atoms at $b\pm \sqrt{b^2+4a}$.}
\end{Example}
Note that the computations of the Cauchy transforms of strongly matricially free convolutions differ 
from those in [12,13], where certain examples of matricially free (but not strongly matricially free) 
convolutions arising in the central limit theorem were studied. 
\noindent\\[10pt]
{\bf Acknowledgement}\\
I would like to thank the anonymous referee for inspiring remarks and suggestions
which were very helpful in the preparation of the revised version of the manuscript.

\end{document}